\mag=\magstep1
\documentstyle{amsppt}
\def\Stab{\operatorname {Stab}}
\def\Hom{\operatorname {Hom}}
\def\dim{\operatorname {dim}}
\def\exp{\operatorname {exp}}
\def\rank{\operatorname {rank}}
\def\<{\langle\langle}
\def\>{\rangle\rangle}
\def\tr{\operatorname {tr}}
\def\ord{\operatorname {ord}}
\def\Im{\operatorname {Im}}
\def\wt{\operatorname {wt}}
\def\Coker{\operatorname {Coker}}
\def\Ker{\operatorname {Ker}}

\def\Spec{\operatorname {Spec}}
\def\Proj{\operatorname {Proj}}
\def\deg{\operatorname {deg}}
\def\Aut{\operatorname {Aut}}
\def\mod{\text{ mod }}
\document
\topmatter
\title
Hodge and Tate conjectures for Hypergeometric sheaves
\endtitle
\author
Tomohide Terasoma
\endauthor
\affil
Department of Mathematical Science,
University of Tokyo,3-8-1, Komaba, Meguro, 
Tokyo, 153, Japan.
\endaffil
\endtopmatter
\heading
\S 0. Introduction
\endheading
In this paper, we study an analog of Hodge and Tate
conjectures for local systems associated to hypergeometric functions.  
Let $f:\Cal X \to S$, 
$g: \Cal Y \to S$ be a smooth proper family of relative dimension $n$ and
$m$ over a smooth variety $S$ over $k=\bold C$.
Let $K$ be an algebraic subfield of $\bold C$.
The higher direct image sheaf $R^if_* K$
and $R^j g_* K$ are variations of $K$-Hodge
structures.  
(For the definition of $K$-Hodge structure, see Section 3.)
A flat proper family $\Cal Z$ of subvarieties in 
$\Cal X \times_S \Cal Y$ of codimension $n$ defines the following
morphism. 
$$
R^if_*K \to R^i(f\times g)_* K \overset{\cap \Cal Z}\to\to
R^{i+2d}(f\times g)_* K \to R^{i}g_* K
$$
This is a morphism of variation of $K$-Hodge structures.
A formal $K$-linear combination a flat families of subvarieties in
$\Cal X \times \Cal Y$ is called a $K$-algebraic correspondence
from $\Cal Y$ to $\Cal X$. By extending $K$-linearly, a 
$K$-algebraic correspondence induces a morphism $R^if_*K \to R^ig_*K$
of variation of $K$-Hodge structures.
\proclaim{Hodge Conjecture for families of varieties}
A morphism 
$R^if_*\bold Q \to R^{i}g_* \bold Q$
of variation of $\bold Q$-Hodge structures over $S$
is induced by a $\bold Q$-algebraic correspondence.
\endproclaim
If $S$ is a point, it is classical Hodge conjecture.
If there exists an action of cyclic group $\mu_m$ on
$\Cal X$ and $\Cal Y$ over $S$, we can define $\chi$-part 
and $\chi '$-part
of the cohomology
group $R^if_*K(\chi )$ and $R^jg_* K(\chi')$ of
$R^i f_*K$ and $R^jg_*K$, where $K = \bold Q (\mu_m)$. Then these
groups are variations of $K$-Hodge structures.  
The $K$-Hodge conjecture
for families of varieties is formulated as follows.
\proclaim{$K$-Hodge Conjecture for families of varieties}
A morphism $R^if_*K(\chi ) \to R^ig_*K(\chi' )$
of variation of $K$-Hodge structures
is induced by a $K$-algebraic correspondence.
\endproclaim
In this paper, we study a local system arising from
the theory of hypergeometric functions due to
Gel'fand-Kapranov-Zelevinski. They are called hypergeometric
sheaves.  Let us recall the definition of hypergeometric sheaves.
Let $n$ and $r$ be natural numbers such that
$r>n+1$ and $\omega_1, \dots,\omega_r$ be elements in 
$\bold Z^n$. For an element $v=(v_1, \dots ,v_n)\in \bold Z^n$,
we define $x^v=x_1^{v_1}\cdots x_n^{v_n}$.
Let $f(a,x) = \sum_{i=1}^r a_ix^{\omega_i}$
be a Laurent polynomial in 
$\bold C[a_1^{\pm}, \dots,a_r^{\pm}, x_1^{\pm}, \dots ,x_n^{\pm}]$
and $U$ be the open set of $(\bold C^\times)^{n+r}$
defined by $U =\{ (x,a) \mid f(a,x) \neq 0\}$.
The covering of $U$ defined by 
$u^m = f(a,x)^{-b_0}x_1^{b_1}\cdots x_n^{b_n}$ is denoted by $\Cal X$.
The group $\mu_d$ acts on $\Cal X$ by the multiplication
on $U$ and
the character $\mu_d \to \bold Q(\mu_d)$ defined by the natural 
inclusion is denoted by $\chi$.  
Then we have a constructible sheaf $\Cal F = R^nf_*K(\chi)$,
where $f : \Cal X \to (\bold C^{\times})^r$
is induced by the second projection $(x,a) \to a$.
We can show that the sheaf $\Cal F \otimes 
K(a_1^{\kappa_1} \cdots a_r^{\kappa_r})$ is defined on a quotient
torus $(\bold C^\times)^{r-n-1}$ of $(\bold C^\times)^r$ for
a suitable choice of $\kappa_1, \dots , \kappa_r \in \bold Q$.
The descent sheaf is called a hypergeometric sheaf.  We can prove that
the hypergeometric sheaf $\Cal F$ is a local system
and has a pure variation of Hodge structure of weight $n$ for an open set 
$U$ in $(\bold C^\times)^{r-n-1}$.
The main theorem of this paper is stated as
follows.
\proclaim{Main Theorem (Theorem 4 (1) in Section 10)}
Let $\Cal F$ and $\Cal F'$ be hypergeometric sheaves 
(for the definition of hypergeometric sheaves see Section 2) on
an open set $U$ in a torus $\bold T$.
If there exists an isomorphism
$\phi : \Cal F \to \Cal F'$ 
of variation of $K$-Hodge structures on $U$,
then there exists an algebraic correspondence $Al$ and a
Hodge cycle of Fermat motif $F$ such that $\phi = Al\cdot F$.
\endproclaim
On the other hand, there is an analogous conjecture, Tate conjecture
for families of varieties.  Let $\bold F_q$ be a finite filed with
$q$-elements and $p$ be its characteristic.  Let $l$ be a prime
number prime to $p$ and $K = \bar\bold Q_l$.  Let $S$ be
a smooth variety over $\bold F_q$ and
$\Cal X \to S$ and $\Cal Y \to S$ be proper smooth varieties over $S$.  Assume
that a cyclic group $\mu_d \subset \bold F_q^{\times}$
acts on $\Cal X$ and $\Cal Y$ over $S$.  Then for a character
$\chi$ and $\chi'$ of $\mu_d$ with the values in $K^\times$,
we can consider the local system of $K$ vector spaces
$R^if_*K(\chi )$ and $R^ig_*(\chi')$ in etale topology.
For a proper flat family $\Cal Z$ of $\Cal X \times_S \Cal Y$,
we can define a homomorphism 
$R^if_*K(\chi ) \to R^ig_*(\chi')$ as an etale sheaf on $S$
in the same way.  Therefore a $K$-algebraic correspondence 
induces a homomorphism $R^if_*K(\chi ) \to R^ig_*(\chi')$ of etale 
sheaves.  Tate conjecture for families of varieties is formulated
as follows.
\proclaim{Tate Conjecture for families of varieties}
A morphism $R^if_*K(\chi ) \to$ 
\linebreak
$R^ig_*K(\chi' )$
of etale sheaves is induced by a $K$-algebraic correspondence.
\endproclaim
Concerning Tate conjecture, we have the following similar
result.
\proclaim{Main Theorem (Theorem 4 (2) in Section 10)}
Let $\Cal F$ and $\Cal F'$ be hypergeometric sheaves 
(for the definition of hypergeometric sheaves see Section 2) on
an open set $U$ in a torus $\bold T$.
If there exists an isomorphism
$\phi : \Cal F \to \Cal F'$ 
of etale sheaves on $U$,
then there exists an algebraic correspondence $Al$ and a
Tate cycle of Fermat motif $F$ such that $\phi = Al\cdot F$.
\endproclaim

To prove the main theorem, we 
recover hypergeometric data from
hypergeometric sheaves.
In this direction, the theory of Mellin transform is useful.
To apply the frame work of Mellin transform, we must overcome
the following points.
\roster
\item
Define Mellin transform as an invariant of local systems.
\item
Find enough informations from Mellin transform.
\item
Construct enough algebraic correspondences to generate the equivalence 
relation arising from Mellin transform.
\endroster
We will discuss these points for $k =\bold C$ and $k=\bold F_q$.
The case $k=\bold C$.
\roster
\item
For a local system on an open set, one can naturally extend to
a perverse sheaf.  They have a structure of cohomological mixed Hodge
complex.  Using this, the Mellin transform is an invariant of
local systems. The argument is a combination of [BBD] and [D].
\item
Hodge numbers provide us enough ample informations to 
detect hypergeometric data up to some equivalence.
This part is quite easy. (See Section 9 Proposition 9.1(1).)
\item
We construct two kinds of algebraic correspondences; one is constant
correspondence and the other is multiplicative correspondence.
In [GKZ], they assume the condition $\sum_{i=1}^r\kappa_i \neq 0$
which is not preserved by constant algebraic correspondences.
We study in more general setting. As a consequence, we study
all together 4 types of constant correspondences.  It is unified
in confluent hypergeometric functions.  We do not take this approach
since they are not algebraic.
\endroster

The case $k=\bold F_q$.
\roster
\item
We should also extend the local system to torus by using perverse
sheaf.  Except the argument for cohomological mixed Hodge complex, it 
is similar to the case $k =\bold C$.
\item
We use the $p$-adic order of the Frobenius on Mellin transform.
It is enough ample to detect the equivalence class for 
hypergeometric data.  In this case,
$p$ power Frobenius action on characters
gives another relation. A typical example is a relation
$g(\chi^p, \psi )= g(\chi, \psi)$.
(see Section 8 for the notation of Gaussian sum.)
Existence of this equality give rise to an
subtle arithmetic relations between $p$-adic value of cohomological 
Mellin transforms. To control these relations, we need Proposition 9.1.(2).
We prove Proposition 9.2 (2) 
in Appendix. By the existence of this equivalence, it is
not natural to assume $\sum_il_i = 0$, which is necessary
for the case $k =\bold C$.
\item
In the case $k=\bold F_q$, two kinds of correspondences;
constant correspondence and multiplicative correspondence is not
sufficient.  We need Frobenius correspondence for hypergeometric data.
\endroster

Let us explain the contents of this paper.  In Section 2, we introduce
hypergeometric data and hypergeometric sheaves.  Hypergeometric
data is a triple 
\linebreak
$D =(R, \{l_i\}_i, \{ \kappa_i\}_i)$
where $R$ is a finitely generated $\bold Z$-module,
$l_i$ is a homomorphism $l_i : R \to \bold Z$ and 
$\kappa_i$ is an element of $I_m(\bar k, K) =
\Hom (\mu_m (\bar k), \mu_m(K))$.  For a hypergeometric data $D$,
we define a constructible sheaf $\Cal G(D)$ 
(resp. $\Cal G(D, \psi )$) on $\bold T(R) =\Spec (k[R])$
if $k= \bold C$ (resp. $k=\bold F_q$).
The sheaves $\Cal G(D)$ and $\Cal G(D, \psi )$ are called 
hypergeometric sheaves.
At the end of Section 2, we prove the perversity and irreducibility
of hypergeometric sheaves.  
The proof is analogous to that in [GKZ].
We review the
definition of variation of $K$-Hodge structures in Section 3.
We introduce a variation of $K$-Hodge structures on the restriction
of hypergeometric sheaves $\Cal G(D)$ to an open set $U$ in 
$\bold T(a)$.  We define a multiplicative equivalence in the set of
hypergeometric data in Section 4.  We show that the resonance condition is
stable under this equivalence relation.  Multiplicative relation
is related to Gauss multiplication formula (See [T].):
$$
\Gamma (ns) =
n^{ns-1}
\frac
{\prod_{i=0}^{n-1}\Gamma (s+\frac{i}{n})}
{\prod_{i=1}^{n-1}\Gamma (\frac{i}{n})}.
$$
This relation produces non-trivial relation between 
hypergeometric sheaves.  In Section 6, these relations are proved
to be induced by an algebraic correspondence.  Before constructing
these algebraic correspondences, we introduce constant equivalence
relation and related algebraic correspondences in Section 5.
The algebraic correspondence for constant equivalence has different
aspects for $k=\bold C$ and $k=\bold F_q$.  In the case
$k = \bold C$, we should discuss algebraic correspondence 
separately according to $\sum_i \kappa_i \neq 0$ or
$\sum_i\kappa_i = 0$.  Here we use the irreducibility
of $\Cal G(D)$ and $\Cal G(D, \psi)$.  In Section 7, we give algebraic
correspondences related to Frobenius action on characters.
Since this equivalence relation does not respect the condition,
$\sum_{i=1}^r l_i =0$, we do not assume this condition if 
$k=\bold F_q$.  In Section 8, we define the cohomological Mellin
transform.  The key proposition, Proposition 9.1 is proved in Section 9
($k =\bold C$) and Appendix ($k = \bold F_q$).  The proof for
the case $k= \bold F_q$ is inspired by the work of [KO] and [A].

Hypergeometric sheaves for the case $k=\bold C$ and $k=\bold F_q$ 
have similar aspects in many points. Nevertheless, we do not work
with the same category.  If $k=\bold C$, we treat hypergeometric sheaves 
with only regular singularities. On the other hand, 
if $k=\bold F_q$, we treat hypergeometric
sheaves with wild ramification.
\demo{Acknowledgement}
The most of this work was done during JAMI project in Johns Hopkins
University.  The author express his thanks for the hospitality.
\enddemo

\heading
\S 1. Notations
\endheading

Let $k$ be the complex number filed $\bold C$ or a finite field $\bold F_q$
of $q$ elements and $p$ be the characteristic of $k$.
Let $V$ be a vector space over $\bold Q$.
A free $\bold Z$ module
$L$ of $V$ is called a lattice if $L \otimes \bold Q = V$.
A finitely generated free $\bold Z$ module is also called a 
lattice. For a lattice $L$ of $V$,
$\bold T (L)$ denotes a torus defined by
$\Spec (k [L])$, where $ k[L]$ is a group ring
over $k$.  
Let $\bold Z_{(p)}$ be the localization of $\bold Z$ at $p$.
For an element $v \in L\otimes \bold Z_{(p)}$, 
the corresponding Kummer character 
on $\bold T(L)$ is denoted by $\kappa (\bold T;v)$ or
$\kappa (v)$ for short.
$K$ denotes the field $\cup_N\bold Q(\mu_N)$
where $\mu_N$ denotes the group of roots of $N$-th power of unity
if $k= \bold C$ and $\bar\bold Q_l$ if $k = \bold F_q$, where
$l$ is a prime number prime to $p$.
The rank 1 local system on $\bold T(L)$
corresponding to the Kummer character $\kappa (\bold T(L),v)$
is denoted by $K(\bold T(L);v)$
or $K(v)$.  The sheaf 
$K(\bold T(L);v)$ is a trivial sheaf of $\bold T(L)$ if
$v \in L$ and for any sublattice $M \subset L$, 
$K(\bold T(L);v)$ has a canonical descent data over
$\bold T (M)$.
For a subvariety $Z$ of $\bold T(L)$, the restriction of 
$K(\bold T(L), v)$ to $Z$ is denoted by $K(Z, v)$.
For a set of variable, for example $x=(x_0, \dots , x_n)$,
$\bold T (x)$ and $\bold A (x)$ denotes a torus and an affine space 
$\Spec k[x_0^{\pm}, \dots , x_n^{\pm}]$ and
$\Spec k[x_0, \dots , x_n]$ whose coordinates
are given by the corresponding variables.

  For a real number $\alpha$, $<\alpha >$ denotes an element
in $\bold R$ such that $<\alpha > \equiv \alpha \quad (\mod \bold Z)$
and $0 \leq \alpha < 1$.  
For fields $E, F$, $I_m(E,F)$ denotes
the set 
\linebreak
$\Hom(\mu_m(E), \mu_m(F))$ and $I(E,F)$ denotes the inductive
limit of $I_m(E,F)$ for all $m$. Then $I(\bold C, \bold C)$
is canonically isomorphic to $\bold Q/\bold Z$ and 
$I(\bar\bold F_q, \bar\bold Q_l)$ is isomorphic to 
$\bold Z_{(p)}/\bold Z$ if we fix a extension $\bar p$ of the ideal $(p)$ 
to the algebraic closure of $\bold Q$ in $\bar\bold Q_l$.

\heading
\S 2 Hypergeometric data and hypergeometric sheaves
\endheading
\heading
\S 2.1 Hypergeometric data
\endheading

In this section, we define hypergeometric sheaves associated to
a data 
\linebreak
$(R, \{ l_i\}_i, \{\kappa_i\}_i)$ called hypergeometric data.
Let $n,r$ be positive integers such that $r-n-1 >0$ and $R$
be a lattice of rank $r-n-1$.  A set $\{ l_i\}_{1,\dots ,r}$
of homomorphism $l_i: R \to \bold Z$ is called primitive if
the induced homomorphism $\bold l=(l_1,\dots , l_r)$
is a primitive injection, i.e. $R$ is identified with a submodule of
$\bold Z^r$ via the homomorphism $\bold l$ whose cokernel is 
torsion free. The homomorphism is called separated if for all $i \neq j$,
there exists no $r \in R$ such that $l_i (r) =1, l_j(r) =-1$
and $l_l(r)=0$ for all $l \neq i,j$.
\demo{Definition}
Let $k$ be $\bold C$ or $\bold F_q$.
The triple $D=(R, \{ l_i\}_{i=1, \dots, j},\{ \kappa_i\}_{i=1, \dots ,r})$
is called a hypergeometric data if 
\roster 
\item
A set of homomorphism $l_i:R \to \bold Z$ $(i=1, \dots , r)$ is 
primitive and separated.
\item
If $k =\bold C$, we assume
$\sum_{i=1}^r l_i = 0$.
\item
$\kappa_i \in I(k, K)$ ($i=1, \dots , r$).
\endroster
\enddemo

In this section we construct a sheaf 
$\Cal G(R,\{ l_i \}_i, \{\kappa_i\}_i) = \Cal G(D)$
or 
\linebreak
$\Cal G(R,\{ l_i \}_i, \{\kappa_i\}_i,\psi) =\Cal G(D, \psi)$
called a hypergeometric
sheaf on a torus $\bold T (R)$.  Let $L$ be the cokernel of
the morphism $\bold l =(l_1, \dots ,l_r): R \to \bold Z^r$.
$$
0 \to R \to \bold Z^r \to L \to 0.
$$
The element in $k[\bold Z^r]$ corresponding to the $i$-th canonical
generator is denoted by $a_i$.
The natural homomorphism $\bold Z^r \to L$ is denoted by $q$.
The image $q(e_v)$ of $e_i = (0, \dots ,
\overset{i}\to{1} , \dots , 0)$ in $L$ is denoted by 
$\omega_i \in L$.  
For an element $v \in L$, the corresponding 
element in $k [L]$ is denoted by $x^v$.
The image $(q\otimes 1)(\kappa )$ of 
$\kappa = (\kappa_i)_i \in \bold Z^r \otimes I(k, K)$
under the natural homomorphism 
$q \otimes 1 : \bold Z^r \otimes I(k, K) \to
L \otimes I(k, K)$ is denoted by $\alpha$.

\heading
\S 2.2 Hypergeometric sheaf for $k=\bold C$.
\endheading

By the condition (2) of hypergeometric data, the homomorphism given
by the summation $\Sigma : \bold Z^r \to \bold Z ;
(u_1, \dots ,u_r) \mapsto \sum_{i=1}^ru_i$
factors through $L$.  The induced map from $L$ to
$\bold Z$ is denoted by $h$.  The kernel of $h$ is denoted by $L_0$.
Then $\bold C [L_0]$ can be identified with a subring of 
$\bold C [ L ]$. 
The element
in $\bold C [\bold Z]$ and $\bold C[L]$, corresponding to 1 and $v \in L$
is denoted by $x_0$ and $x^v$ respectively.
\proclaim{Proposition 2.1}
Let $v \in L$ such that $h(v) >0$ and $a_1, \dots a_r$
be elements in $\bold C^\times$.  Then for a Laurent polynomial
$f = \sum_{i=1}^r a_i x^{\omega_i}$, $(f)^{h(v)}x^{-v}$
is an element in $\bold C[L]$.
\endproclaim
\demo{Proof}
Since $h(v)>0$, it is an element in $\bold C[L]$. The subring
$\bold C [L_0]$ in $\bold C [L]$ is characterized as follows.
For a $\bold C$-algebra homomorphism $\mu : \bold C[\bold Z] \to \bold C$,
the map $x^w \mapsto \mu (x_0)^{h(w)}x^w$ defines a ring endomorphism
$T_{\mu}: \bold C [ L] \to \bold C [ L]$
of $\bold C [ L]$.  The subring $\bold C [L_0]$ is characterized
by
$$
\{ f \in \bold C [ L] \mid
T_{\mu}(f) = f \text{ for all }\mu : \bold C [\bold Z ]\to \bold C \}.
$$
Since
$$
\align
T_{\mu}(f^{h(v)}x^{-v}) = & (\sum_i a_iT_{\mu}(x^{\omega_i}))^{h(v)}
\mu (x_0)^{-h(v)} x^{-v} \\
= & (\sum_i a_i\mu (x_0)x^{\omega_i})^{h(v)}
\mu (x_0)^{-h(v)} x^{-v} \\
= & f^{h(v)}x^{-v}, \\
\endalign
$$
we obtain the proposition.
\enddemo
\proclaim{Corollary}
Let $f$ be as in Proposition 2.1. The ideal $(f)$ in $\bold C [L]$
is defined over $\bold C[L]$, i.e. if we put 
$I_f = (f) \cup \bold C [L_0]$, then $I_f\bold C [L] = (f)$.
\endproclaim

Now we consider the relative situation.  Let $a$ be a system of
variables $a = (a_1, \dots , a_r)$.  Let $f$ be a polynomial
in $\bold C [a_i^{\pm}] \otimes \bold C [ L]$ defined by
$$
f = \sum_{i=1}^r a_i x^{\omega_i}.
$$
\heading
Case $h(\alpha ) \neq 0$
\endheading

The subvariety of $\bold T (a) \times \bold T( L)$
defined by $\{ f = 1\}$ is denote by $\Cal X$ and
the composite $\Cal X \to \bold T(a)$ of the natural inclusion 
$\Cal X \to \bold T (a)\times \bold T ( L)$ and
the first projection 
$\bold T (a) \times \bold T (L) \to \bold T (a)$
is denoted by $\varphi$.
The element $x^\alpha$ defines a rank 1 local systems
$K(\bold T( L);x^\alpha )$ and
$K(\Cal X ;x^\alpha )$
of $K$-vector spaces
over $\bold T ( L)$ and $\Cal X$.
The $i$-th higher direct image sheaf
of $K(\Cal X ;x^\alpha )$ by $\varphi$ is denoted by 
$R^i\varphi_* K(\Cal X ; x^\alpha )$. 
\heading
Case $h(\alpha )=0$
\endheading

By the Corollary to Proposition 2.1, the variety 
$\Cal X' =\{ (x,a) \in \bold T(L) \times \bold T (a) \mid f(a,x)=0\}$ 
is defined over $\bold T(L_0)$, i.e. there is unique subvariety
of $\bold T(L_0)\times\bold T(a)$ such that $\Cal X'$ is
isomorphic to the pull back of $\Cal X$.
Under the condition $h(\alpha )=0$,
we show that the rank 1 local system 
$K (\prod_{i=1}^r (x^{\omega_i})^{\kappa_i})$ 
corresponding to the Kummer character
$\kappa (\bold T(L) \times \bold T(a);
\prod_{i=1}^r (x^{\omega_i})^{\kappa_i})$ is naturally descent to 
$\bold T(L_0) \times \bold T(a)$ via the the morphism
$\bold T(L) \times \bold T(a) \to \bold T(L_0) \times \bold T(a)$.
To prove this, it is enough to show that the element
$\prod_{i=1}^r (x^{\omega_i})^{\kappa_i}$ is invariant
under the automorphism $T_{\mu}$ for
any algebra homomorphism $\mu : \bold C [\bold Z] \to \bold C$.
Since 
$$
\align
T_{\mu}(\prod_{i=1}^r (x^{\omega_i})^{\kappa_i}) = &
(\prod_{i=1}^r (\mu (x_0)x^{\omega_i})^{\kappa_i})  \\
= & (\prod_{i=1}^r (x^{\omega_i})^{\kappa_i}),  \\
\endalign
$$
we have the claim. The restriction of the descent sheaf
to the variety $\Cal X$ is denoted by $K(\Cal X ;x^{\alpha})$.
Let $\varphi $ be the composite of the closed immersion
$\Cal X  \to \bold T (L_0) \times \bold T(a)$
and the second projection $\bold T(L_0) \times \bold T(a) \to \bold T(a)$.
The $i$-th higher direct image of $K(\Cal X ;x^{\alpha})$ is denoted
by $R^i\varphi_{*}K(\Cal X ;x^{\alpha})$.

The element 
$a^\kappa = (a_i^{\kappa_i})_i \in \bold Z^r \otimes I(\bold C, K)$
defines a rank 1 local system 
of $K$-vector spaces $K(\bold T(a) ;a^{\kappa} )$.
We descend the the sheaf 
$R^i\varphi_*K(\Cal X ;x^\alpha ) \otimes K(\bold T(a) ;a^{\kappa} )$
to $\bold T(R)$ via the natural morphism $\bold T(a) \to \bold T(R)$
arising from the homomorphism $\bold l : \bold R \to \bold Z^r$.
Since we have the natural exact sequence
$$
\CD
0 @>>> R @>>> \bold Z^r \oplus \bold Z^r @>>> 
\bold Z^r \oplus L @>>> 0 \\
& & r &\mapsto &(r, -r) \\
& & & & (a, b) &\mapsto  & (a+b, q(b)) \\
\endCD
$$
The fiber product $\bold T (a) \times_{\bold T(R)} \bold T (a)$
is identified with $\bold T (a) \times \bold T ( L)$ and
the following diagram commutes.
$$
\CD
& \bold T (a) \times_{\bold T(R)} \bold T (a) \\
 pr_1 \swarrow & & \searrow pr_2 \\
\bold T (a) &  @V{\iota}V{\simeq}V \bold T (a) \\
 pr_1 \nwarrow & & \nearrow  m \\
& \bold T(a) \times \bold T( L) \\
\endCD
$$
Here $m$ is given by the homomorphism of the modules
$$
\bold Z^r \to \bold Z^r \oplus L ; u \mapsto (u, q(u)).
$$
\proclaim{Proposition 2.2}
Let $\tilde\Cal G = \tilde\Cal G(D) = 
\bold R^n\varphi_*K(\Cal X ;x^\alpha )\otimes K(\bold T(a) ;a^{\kappa} )$
or 
\linebreak
$\bold R^{n-1}\varphi_*K(\Cal X ;x^\alpha )\otimes K(\bold T(a) ;a^{\kappa} )$
according to $\sum_{i=1}^r \kappa_i \neq 0$ or
$\sum_{i=1}^r \kappa_i = 0$
Then there exists a natural isomorphism $\phi$
$$
\phi :m^* (\tilde\Cal G) \simeq  pr_1^*(\tilde\Cal G)
$$
which satisfies the cocycle condition for descent data via 
the isomorphism
$\iota$.
\endproclaim
\demo{Proof}
First we consider base change diagrams of $\Cal X$:
$$
\CD
\Cal X_1 @>>> \Cal X @<<< \Cal X_2 \\
@VVV @VVV @VVV \\
\bold T(a) \times \bold T ( L) @>>{m}> \bold T(a) @<<{pr_1}< 
\bold T(a) \times \bold T ( L) \\
\endCD
$$
Here $pr_1^*,m^*$ are the ring homomorphisms
$$
\bold C [a_i^{\pm}] \to \bold C[a_i^{\pm},y^v]_{v\in L}
$$
given by $pr_1^*(a_i) =a_i$ and $m^*(a_i)=a_iy^{\omega_i}$
with a new set of variable $y^v$ $(v \in L)$.  If 
$\sum_{i=1}^r \kappa_i \neq 0$ (resp. $\sum_{i=1}^r\kappa_i =0$),
the defining equation of $\Cal X_1$ and $\Cal X_2$ is given by
$$
\align
\Cal X_1 :  \sum_{i=1}^r a_ix^{\omega_i} = 1, & \quad(\text{ resp. }
\Cal X_1 :  \sum_{i=1}^r a_ix^{\omega_i} = 0 ) \\
\Cal X_2 :  \sum_{i=1}^r a_iy^{\omega_i}x^{\omega_i} = 1, & 
\quad(\text{ resp. }
\Cal X_2 :  \sum_{i=1}^r a_i
y^{\omega_i-\omega_1}x^{\omega_i} = 0 ). \\
\endalign
$$
Therefore the automorphisms of 
$\bold C [a_i^{\pm}, x^v, y^v ]_{v\in  L}$ and
$\bold C [a_i^{\pm}, x^{v'}, y^v ]_{v' \in  L_0, v \in L}$
given by 
$$
\cases
\phi^* : a _i \mapsto a_i, x^v \mapsto x^v y^v , y^v \mapsto y^v
 & (\sum_{i=1}^r \kappa_i \neq 0) \\
{\phi'}^* : a _i \mapsto a_i, x^{v'} \mapsto x^{v'} y^{v'} , y^v \mapsto y^v
& (\sum_{i=1}^r \kappa_i = 0) \\
\endcases
$$
induce isomorphisms $\Cal X_2 \to \Cal X_1$.  The base
change of Kummer character associated to 
$a_1^{\kappa_1}\cdots a_r^{\kappa_r}x^{\alpha}
\in (\bold Z^r \oplus  L)\otimes I(\bold C, K)$ on $\Cal X$
to $\Cal X_1$ and $\Cal X_2$ are given by
$a_1^{\kappa_1} \cdots a_r^{\kappa_r} x^{\alpha}$
and
$$
\align
m^*(a_1^{\kappa_1}\cdots a_r^{\kappa_r})x^{\alpha} =
& (\prod_{i=1}^r a_i^{\kappa_i}y^{\kappa_i\omega_i})x^{\alpha} \\
= &  (\prod_{i=1}^r a_i^{\kappa_i})y^{\alpha}x^{\alpha} \\
\endalign
$$
respectively.  This two Kummer character on $\Cal X_1$ and $\Cal X_2$
corresponds to each other via the isomorphism $\phi$ and $\phi '$ 
constructed as before.
\enddemo
\demo{Definition} The sheaf obtained by descending $\tilde \Cal G (D)$
to $\bold T(R)$ is denoted by $\Cal G(D)$.
The sheaf $\Cal G (D)$ is called hypergeometric sheaf for a hypergeometric
data $D= (R, \{l_i\}_i, \{\kappa_i\}_i)$.
Note that, $\tilde\Cal G(D)$ is the equal to 
$R^n\varphi_*
K(\Cal X;\prod_{i=1}^r(a_ix^{\omega_i})^{\kappa_i})$ or
$R^{n-1}\varphi_*
K(\Cal X;\prod_{i=1}^r(a_ix^{\omega_i})^{\kappa_i})$ according to
$\sum_{i=1}^r\kappa_i \neq 0$ or $\sum_{i=1}^r \kappa_i=0$.
\enddemo
We compare the sheaf $\Cal G(D)$ and the sheaf associated to
Gel'fand-Kapranov-Zelevinski (=GKZ for short) hypergeometric functions
by choosing a suitable base of $ L$ under the assumption
$\sum_{i=1}^r \kappa_i \neq 0$.
Let $s: \bold Z \to L$ be a section of $h$ and take a base
of $L_0 \simeq \oplus_{i=1}^n\bold Z v_i$.  Then the ring $\bold C [L_0]$
and $\bold C [ L]$ can be identified with Laurent polynomial 
rings $\bold C [x_1^{\pm},\cdots , x_n^{\pm}]$ and
$\bold C [x_0^{\pm},x_1^{\pm},\cdots , x_n^{\pm}]$, where
$x_i$ corresponds the element $v_i$ of $L$ and $x_0$ corresponds to
the image of $1\in \bold Z$ under the section $s$.
Using the base of $ L$, the element $\omega_i$
is expressed as $(1, \bar\omega_i)$ with $\bar\omega_i \in \bold Z^n$.
Then by the definition of $L$, it is generated by 
$\omega_i$.
We introduce a system of variables $\bar x = (x_1, \dots , x_n)$ and
use the common notation for multiple index: 
$\bar x^{v} = x_1^{v_1}\cdots x_n^{v_n}$ for $v \in \bold Z^n$.
Then $x^{\omega_i}=x_0\bar x^{\bar\omega_i}$ and the defining
equation of $\Cal X$ is $x_0f_0(a_1, \dots ,a_r,x_1,\dots x_n) = 1$,
where $f_0=\sum_{i=1}^r a_i \bar x^{\bar\omega_i}$. If we write 
$$
\alpha =(\alpha_0, \alpha_1, \dots ,\alpha_n) 
\in L\otimes I(\bold C, K) 
\simeq \bold Z^{n+1} \otimes I (\bold C, K),
$$
then the corresponding Kummer character on $\Cal X$ is
$x_0^{\alpha_0}x_1^{\alpha_1}\cdots x_n^{\alpha_n} =$
\linebreak
$f_0^{-\alpha_0}x_1^{\alpha_1}\cdots x_n^{\alpha_n}$,
where $\alpha_0 = h(\alpha ) =\sum_{i=1}^r \kappa_i$
Let $\Cal U$ be the open set of $\bold T (x) \times \bold T(a)$
defined by $\Cal U =\{ (x,a) \in \bold T (x) 
\times \bold T (a) \mid
f_0(a,x) \neq 0 \}$,
where $x =(x_1, \dots x_n)$.  
The rank 1 local system corresponding to the Kummer character
$f_0^{-\alpha_0}x_1^{\alpha_1}\cdots x_n^{\alpha_n}$
is denoted by $K_{\Cal U}(x^\alpha)$. Then $R^i\varphi_* K(x^\alpha )$
is identified with $R^i{\bar\varphi}_* K_{\Cal U}(x^\alpha )$, where
$\bar\varphi :\Cal U \to \bold T (a)$ is the composite of the open immersion
to $\bold T(x) \times \bold T(a)$ and the second projection.
This is nothing but the local system defined by 
Gel'fand-Kapranov-Zelevinski.  Note that is it is independent of
the choice of the base $L$ given as above.

\heading
\S 2.3 Hypergeometric sheaf for $k=\bold F_q$
\endheading
In this case,
we use an extra parameter $t$.
Let $\Cal X$ be the subvariety of $\bold T(x) \times \bold A (t)
\times \bold T(a)$ defined by;
$$
\Cal X = \{ (x,t ,a) \in \bold T(x) \times \bold A (t)
\times \bold T(a) \mid
f(a, x)=t \},
$$
where $f(a,x) = \sum_{i=1}^ra_i x^{\omega_i}$.
Let $\varphi : \Cal X \to \bold T(a)$ be the composite of
the natural immersion $\Cal X \to \bold T (x) \times \bold A (t )
\times \bold T (a)$ and the third projection 
$\bold T (x) \times \bold A (t )
\times \bold T (a) \to \bold T(a)$.
Let $\psi$ be a non-trivial additive character 
$\psi : \bold F_p \to \bar\bold Q_l$.  
Let $\pi :\bold A(\tau) \to \bold A(t)$ be the 
Artin-Shreier covering defined by $t= \tau^p - \tau$.  The character
$\psi$ of the covering transformation group $\bold F_p$ defines
a rank 1 local system $\Cal L_{\psi}$ of $\bold A(t)$.
This local system is called the Lang sheaf.
The pull back of $\Cal L_{\psi}$ by the projection to the
$t$-line is denoted by $\Cal L_{\psi}(t)$. 
We can prove the following proposition in the same way.
\proclaim{Proposition 2.3}
The sheaf 
$R^n\varphi_*(K(\Cal X ;x^{\alpha}a^{\kappa})\otimes \Cal L_{\psi}(t ))=
R^n\varphi_*(K(\Cal X ;x^\alpha )\otimes \Cal L_{\psi}(t ))
\otimes K(\bold T(a) ;a^\kappa)$
on $\bold T(a)$
can be canonically descent to a constructible 
sheaf $\Cal G(D,\psi )$ on $\bold T(R)$.
\endproclaim
\demo{Definition}
The sheaf $\Cal G(D, \psi )$ is called the hypergeometric sheaf for
a hypergeometric data $D$ if $k =\bold F_q$.
\enddemo

\heading
\S 2.4 Properties of hypergeometric sheaves
\endheading

Now we recall several properties 
on $R^i\varphi_*K(\Cal X ;x^\alpha )$.  Let $\Delta$
be the convex hull of 
$\omega_i$ and $0$ and $C(\Delta)$ be the convex cone
generated by $\Delta$ in
$L \otimes \bold R$.
For a codimension 1 face $\sigma$ of $C( \Delta )$, $H_{\sigma}$
denotes the linear hull of $\sigma$ and $h_{\sigma}$ denotes a
primitive linear form defining $H_{\sigma}$.  
An element $\alpha \in L\otimes I(k,K)$ 
is called non-resonant if
$h_{\sigma}(\alpha )$ is
not zero.
(Note that if $\Delta$ contains 0 as interior, there is no codimension
1 face in $C(\Delta )$).
For a face $\sigma$ of $\Delta$, we define $f_{\sigma}$ by
$$
f_{\sigma} = \sum_{\omega_i \in \sigma }a_i x^{\omega_i}.
$$
Let $a^{(0)}\in \bold T(a)$.
The polynomial $f^{(0)}(x)=f(a^{(0)},x)$ 
is said to be non-degenerate 
if the variety
$\{ f_{\sigma}^{(0)} = 0\}$ is smooth in $\bold T[ L]$ for all the
faces $\sigma$ of $\Delta$. 
A point $a^{(0)}$ is called non-degenerated if $f^{(0)}$ is 
non-degenerated.
The open set consisting
the non-degenerate point in $\bold T(a)$ is denoted by $U_a$.  
It is easy to see that for a point $r^{(0)} \in \bold T (R)$, all the
points in $\bold l^{-1}(r^{(0)})$ is non-degenerated if
a point $a^{(0)}$ in $\bold l^{-1}(r^{(0)})$ is non-degenerate.
A point $r^{(0)}$ in $\bold T(R)$ is called non-degenerate if there
exists a non-degenerate point $a^{(0)}$ in $\bold l^{-1}(r^{(0)})$.
The open set of $\bold T(R)$ consisting of non-degenerate point
is denoted by $U_R$.

For $k=\bold C$ or $\bold F$, we have the following proposition.  The
related statement for $\Cal D$-module, it is proved in [GKZ].
\proclaim{Proposition 2.4}
\roster
\item
Let $k=\bold C$.
If $\alpha$ is non-resonant, then the
restriction of $R^n\varphi_*K(\Cal X ;x^{\alpha} )$ (resp.
$R^{n-1} \varphi_*K(\Cal X ;x^{\alpha} )$) to the open set $U_a$ is a 
local system of $K$ vector spaces if $\sum_{i=1}^r\kappa_i \neq 0$
(resp. $\sum_{i=1}^r\kappa_i =0$).
\item
Let $k=\bold F_q$. If $\alpha$ is non-resonant, then the
restriction of $R^n\varphi_*K(\Cal X ;x^{\alpha} )\otimes \Cal L_{\psi}(t )$ 
to the open set $U_a$ is a local system of $K$ vector spaces.
\endroster
\endproclaim
\proclaim{Proposition 2.5}
Suppose that $\alpha$ is non-resonant. 
\roster
\item
Let $k=\bold C$. The natural homomorphism 
$$
R^i\varphi_! K(\Cal X ;x^{\alpha} ) \to R^i \varphi_* K(\Cal X ;x^{\alpha} )
$$
is an isomorphism. As a consequence,
$R^i\varphi_*K(\Cal X;x^{\alpha} )$ is 0 unless $i=n$ 
(resp. $i=n-1$) 
if $\sum_{i=1}^r\kappa_i \neq 0$
(resp. $\sum_{i=1}^r\kappa_i =0$).
\item
Let $k=\bold F_q$. The natural homomorphism 
$$
R^i\varphi_! K(\Cal X ; x^{\alpha} )\otimes \Cal L_{\psi}(t ) 
\to R^i \varphi_* K(\Cal X ;x^{\alpha} )\otimes \Cal L_{\psi}(t )
$$
is an isomorphism. As a consequence,
$R^i\varphi_*K(\Cal X ;x^{\alpha})\otimes \Cal L_{\psi}(t)$
is $0$ unless $i=n$.
\item
Let $k=\bold C$.
The constructible sheaf $R^n\varphi_*K(\Cal X ;x^{\alpha} )$
\linebreak
(resp. $R^{n-1}\varphi_*K(\Cal X ;x^{\alpha} )$) is an irreducible
perverse sheaf. As a consequence, the restriction of 
$R^n\varphi_*K(\Cal X ;x^{\alpha} )$ (resp.
$R^{n-1}\varphi_*K(\Cal X ;x^{\alpha} )$)
to $U_a$ is an irreducible local system of K-vector spaces.
\item
Let $k=\bold F_q$.
The constructible sheaf 
$R^n\varphi_*K(\Cal X ;x^{\alpha} )\otimes\Cal L_{\psi} (t )$
is an irreducible perverse sheaf. As a consequence, the restriction of 
$R^n\varphi_*K(\Cal X ;x^{\alpha} )\otimes \Cal L_{\psi}(t)$
to $U_a$ is an irreducible local system of K-vector spaces.
\endroster
\endproclaim
\demo{Proof}
We use a method of Fourier transform.
First we prove the following Lemma
\enddemo
First we prove the propositions for $k= \bold F_q$.
Let $y^v$ ($v \in L$) and $b=(b_1, \dots ,b_r)$ be systems of variables.
Let $j:\bold T(y) \to \bold A(b)$ be the morphism defined by
$q : \bold N^r \to L$. 
\proclaim{Lemma 2.6}
Under the resonance condition, 
$R^i j_! K(\bold T(y); y^{\alpha})\to 
R^i j_* K(\bold T(y) ;y^{\alpha})$ is an isomorphism.
\endproclaim
We use toric geometry.
Let $\Delta_0$ be the image of the semi-group homomorphism
$\bold N^r \to L$. Then $A_{\Delta_0}=\Spec k[\Delta_0]$ is the closure of
the image of $j$. Let $\Delta$ be the intersection of
the convex hull $\Delta_{\bold R}$
of $q(e_i)=\omega_i$ in $L\otimes \bold R$ and $L$. 
Then $A_{\Delta} = \Spec k[\Delta ]$ is the normalization of $A_{\Delta_0}$.
Let $\tilde F$ be a regular triangulation of the dual fan 
$F=\{C_b \mid b \in \Delta_{\bold R}\}$ 
of $\Delta_{\bold R}$, where $C_b =\{ x \in L^*\otimes \bold R \mid
(x, b'-b) \geq 0 \text{ for all } b'\in \Delta_{\bold R}\}$.  
Then the toric variety
$X(\tilde F)$ is a smooth variety and the natural morphism 
$p :X(\tilde F) \to A_{\Delta}$ is proper.  
Therefore the composite $\pi : X(\tilde F) \to \bold A(b)$ is proper.
Let 
$\tilde j : \bold T(y) \to X(\tilde F)$ be the natural open immersion
and $D$ be the complement of $\bold T(y)$ in $X(\tilde F)$.
$$
\bold T(y) \overset{\tilde j}\to\to X(\tilde F)
\overset{p}\to\to A_{\Delta} \overset{\pi}\to\to
\bold A (b)
$$
To prove the lemma
it is enough to prove 
$\bold R p_!(C) = 0$, where $C$ is the mapping cone;
$C(\bold R \tilde j_! K(\bold T(y) ;y^{\alpha}) 
\to \bold R \tilde j_* K(\bold T(y) ;y^{\alpha}))$.
Therefore, it is enough to prove that
$$
H^i_c(p^{-1}(x), R^k\tilde j_*K(\bold T(y) ;y^\alpha) \mid_{D \cap p^{-1}(x)})=0
\tag{1}
$$
for all $i,k$ and $x \in A_{\Delta_0}$.
  Let $D= \coprod_{\beta \in \tilde A}S_{\beta}$ be
the decomposition of $D$ by toric strata, where $\tilde A$
denotes a set of positive dimensional cone in $\tilde F$.  By the long exact
sequence for cohomologies with compact support, to show (1), 
it is enough to prove that
$$
H^i_c(p^{-1}(x)\cap S_{\beta}, 
R^k\tilde j_*K(\bold T(y) ;y^\alpha) \mid_{S_{\beta} \cap p^{-1}})=0.
\tag{2}
$$
Let $L_{\alpha}^*$ be the intersection of $\bold R$-vector space
generated by $\alpha$ and the lattice $L^*$.  The element
$\alpha \in L\otimes I( k, K)$ defines a homomorphism $\alpha_{\beta}:$
$L_{\alpha}^* \to I( k, K)$.  It is easy to see that  
$$
R^i\tilde j_*K(\bold T(y);y^{\alpha})\mid_{S_{\beta}}=
\cases
 (\wedge^i K^{\dim (\beta )}(-1) )\otimes K(S_{\beta} ;y^{\alpha})
& \qquad \text{ if } \alpha_{\beta}= 0 \\
 0 & \qquad \text{ if } \alpha_{\beta} \neq 0. \\
\endcases
$$
The restriction of $p$ to $S_{\beta} \to A_{\Delta}$ corresponds to
the primitive embedding $L_{\beta '} \to L_{\beta}$, where
$\beta '$ is the minimal cone in $F$ containing $\beta$ and
$L_{\beta}$ and $L_{\beta'}$ is a primitive lattice of $L$
annihilated by all the elements in $\beta$ and $\beta'$ respectively.
If $R^i\tilde j_*K(\bold T(y) ;y^{\alpha})\mid_{S_{\beta}}$
does not vanish,  $\alpha \in L_{\beta}\otimes I( k, K)$.  
On the other hand, by the
resonance condition, $\alpha \notin L_{\beta'}\otimes I( k, K)$.  Therefore, 
$K(S_{\beta} ;y^{\alpha})$ has not trivial monodromy along fibers.
Therefore the equality (2) holds.
This completes the proof of the lemma.
\demo{Proof of Proposition 2.5}
Let $k: \bold A(a) \times \bold T(y) \to \bold A(a) \times \bold A(b)$
be the product of the inclusion given in Lemma 2.6.
Let $\Cal F= \bold R j_{!} K(\bold T(y) ;y^\alpha) = 
\bold R j_{*} K(\bold T(y) ;y^\alpha)$.  Let $\bold P(b)$, $H_{0}(b)$ and $k$ 
be the projective space compactifying $\bold A(b)$,
$\bold P(b)-\bold A(b)$ and the open immersion
$k : \bold A(a) \times \bold A(b) \to \bold A(a) \times \bold P(b)$.
The sheaf $K(\bold T(y) ;y^{\alpha})$ has trivial monodromy on some tame
covering $\bold P(\tilde b)$ of $\bold P(b)$.
Since $\Cal L_{\psi}(\sum_{i=1}^ra_ib_i)$ 
ramifies wildly along $H_{0}$,
the natural homomorphism
$$
\bold R k_!(pr_2^* \Cal F\otimes \Cal L_{\psi}(\sum_{i=1}^ra_ib_i)) \to
\bold R k_*(pr_2^* \Cal F\otimes \Cal L_{\psi}(\sum_{i=1}^ra_ib_i))
$$
is a quasi-isomorphism.
On the other hand, 
$R^i \varphi_!K(\Cal X ;x^\alpha)\otimes \Cal L_{\psi}(t)$ and
\linebreak
$R^i \varphi_*K(\Cal X ;x^\alpha)\otimes \Cal L_{\psi}(t)$ are isomorphic to 
$\bold R {pr_1}_! (\bold R k_!(pr_2^* K(y^\alpha)
\otimes \Cal L_{\psi}(\sum_{i=1}^ra_ib_i)))$ and
$\bold R {pr_1}_* (\bold R k_*(pr_2^* K(y^\alpha)
\otimes \Cal L_{\psi}(\sum_{i=1}^ra_ib_i)))$
respectively, where $pr_1 : \bold A(a) \times \bold P(b) \to \bold A(a)$.
Since $pr_1$ is proper, we get the proposition.
\enddemo

We get Proposition 2.5 for $k=\bold C$ by reduction
mod $p$.
Let $D = (R, \{l_i\}_i, \{ \kappa_i \}_i)$ be a hypergeometric data for
$k=\bold C$. 
By choosing an embedding $\bold Q(\mu_d) \to \bar\bold Q_l$,
the sheaf 
$R^i\varphi_* K(\Cal X ;\prod_{i=0}^nx^\alpha_i) \otimes
\bar\bold Q_l$ is considered as an etale cohomology.
Let $d$ be the common denominator of $\kappa_i$ 
($i= 1, \dots , r$), $p$ be a rational prime prime to $d$
and $\frak p$ be a prime of $\bold Q(\mu_d)$ over $p$.
The residue field $\kappa (\frak p)$ at $\frak p$ is a finite filed
$\bold F_q$.
Since the variety $\Cal X$ is defined over $\bold Q(\mu_d)$,
we can consider the reduction of 
$R^i\varphi_* K(\Cal X ;\prod_{i=0}^nx^\alpha_i) \otimes
\bar\bold Q_l (\text{ mod }\frak p)$ at 
$\frak p$ using obvious model over $\bold Z[\mu_d]$.
On the other hand, the data $D$ is also a hypergeometric
data for $k=\bold F_q$.  The variety defining
hypergeometric sheaf over $\bold F_q$ is denoted by $\Cal X^{(p)}$ 
and the natural morphism $\Cal X^{(p)} \to \bold T(a)$ is denoted
by $\varphi^{(p)}$ for $k=\bold F_q$.
We compare the sheaf 
$R^i\varphi_* K(\Cal X ;\prod_{i=0}^nx^\alpha_i) \otimes
\bar\bold Q_l (\text{ mod }\frak p)$  (resp.
$R^i\varphi_! K(\Cal X ;\prod_{i=0}^nx^\alpha_i) \otimes
\bar\bold Q_l (\text{ mod }\frak p)$) and
$R^i\varphi^{(p)}_* K(\Cal X^{(p)} ;
\prod_{i=0}^nx^\alpha_i)\otimes\Cal L_{\psi}(t )$.
(resp. $R^i\varphi^{(p)}_! K(\Cal X^{(p)} ;
\prod_{i=0}^nx^\alpha_i)\otimes\Cal L_{\psi}(t)$ .)
The following proposition completes the proof of Proposition 2.5
for $k= \bold C$.
\proclaim{Proposition 2.7}
We have the following isomorphism
$$
\align
& R^i\varphi^{(p)}_* K(\Cal X^{(p)};\prod_{i=0}^nx^\alpha_i)(\psi ) \\
& \simeq 
\cases
R^{i-1}\varphi_* K(\Cal X ;\prod_{i=0}^nx^\alpha_i) \otimes
\bar\bold Q_l (\text{ mod } \frak p) \otimes g(\alpha_0, \psi )
& \qquad (\text{ if }\alpha_0 \neq 0)  \\
R^{i-2}\varphi_* K(\Cal X ;\prod_{i=0}^nx^\alpha_i) \otimes
\bar\bold Q_l(-1) (\text{ mod } \frak p)
& \qquad (\text{ if }\alpha_0 = 0)  \\
\endcases  \\
\endalign
$$
$$
\align
& R^i\varphi^{(p)}_! K(\Cal X^{(p)};\prod_{i=0}^nx^\alpha_i)(\psi ) \\
& \simeq 
\cases
R^{i-1}\varphi_! K(\Cal X ;\prod_{i=0}^nx^\alpha_i) \otimes
\bar\bold Q_l (\text{ mod } \frak p) \otimes g(\alpha_0, \psi )
& \qquad (\text{ if }\alpha_0 \neq 0)  \\
R^{i-2}\varphi_! K(\Cal X ;\prod_{i=0}^nx^\alpha_i) \otimes
\bar\bold Q_l(-1) (\text{ mod } \frak p)
& \qquad (\text{ if }\alpha_0 = 0),  \\
\endcases  \\
\endalign
$$
where $g(\alpha_0, \psi)$ is a $Gal (\bar\bold F_q/\bold F_q)$-module
whose action of Frobenius is the multiplication by the Gaussian sum
$$
g(\alpha_0, \psi) = \sum_{x \in \bold F_q^{\times}}\chi_{\alpha_0}(x)
\psi (x),
$$
where $\chi_\alpha$ is a character of $\bold F_q$ defined by
$\chi_\alpha (x)( \text{ mod }\frak p) = x^\frac{q-1}{m}$.

\endproclaim
\demo{Proof}
We take a base of $L\simeq \oplus_{i=0}^n \bold Z$ 
such that the $0$-th projection
coincides with $h$. Under this choice of base, the ring
$k[L]$ is identified with $k[x_0^{\pm}, \dots ,x_n^{\pm}]$. Let
$x$ denotes the system of variable $x= (x_1, \dots ,x_n)$.
Then the image of $e_i$ under the morphism $q : \bold Z_r \to L$
is denoted by 
$(1, \bar\omega_i)$.  We define 
$f_0 = \sum_{i=1}^r a_ix^{\bar\omega_i}$.
Let $\Cal X$ and $\Cal X^{(p)}$ be the variety over
$\bold Q(\mu_d)$, $\bold F_q$ and $W$ defined
by
$$
\align
\Cal X & = \{ (x_0, x, a) \in \bold T(x_0) \times \bold T(x)
\times \bold T(a) \mid
x_0f_0(a,x) =1 \} \\
\Cal X^{(p)} & = \{ (x_0, x, a, t) \in \bold T(x_0) \times \bold T(x)
\times \bold T(a) \times \bold A(t) \mid
x_0f_0(a,x) =t \} \\
W  & = \{ (x_0, y, t) \in \bold T (x_0) \times \bold A(y) \times
\bold A(t) \mid x_0y =t \}. \\
\endalign
$$
Consider the following fiber product;
$$
\CD
\Cal X^{(p)} @>{F'}>> W \\
@V{\pi '}VV @VV{\pi}V \\
\bold T(a) \times \bold T(a) @>>{F}> \bold A(y)  \\
(a,x) &\mapsto & f_0(a,x) \\
\endCD
$$
The pull back  $pr_3^*\Cal L_{\psi}$ 
of the Lang sheaf $\Cal L_{\psi}$ on $\bold A(t)$
under the morphism $W \to \bold A(t)$ is also denoted by $\Cal L_{\psi}$
for short. Let $K(W ;x_0^{\alpha_0})$ be the 
rank 1 local system corresponding to
the Kummer character $x_0^{\alpha_0}$ on $W$. Then 
$K(W ;x_0^{\alpha_0}) \otimes \Cal L_{\psi}$ is a local system on $W$.
Since the morphism $F$ is smooth, we have
$$
\align
F^*\bold R \pi_*( K(W ;x_0^{\alpha_0}) \otimes \Cal L_{\psi}) & \simeq
\bold R \pi'_*(F')^*(K(W; x_0^{\alpha_0}) \otimes \Cal L_{\psi}) \\
F^*\bold R \pi_!( K(W;x_0^{\alpha_0}) \otimes \Cal L_{\psi}) & \simeq
\bold R \pi'_!(F')^*(K(W;x_0^{\alpha_0}) \otimes \Cal L_{\psi}) \\
\endalign
$$
The sheaf
$R^i \pi_*( K(W;x_0^{\alpha_0}) \otimes \Cal L_{\psi})$ and
$R^i \pi_!( K(W;x_0^{\alpha_0}) \otimes \Cal L_{\psi})$
are computed as follows.
If $\alpha_0$ is not trivial,
$$
\align
R^i \pi_*( K(W;x_0^{\alpha_0}) \otimes \Cal L_{\psi})
& =
\cases
\tilde j_*K(y^{-\alpha_0}) \otimes g(x_0^{\alpha_0}, \psi ) 
& \qquad (\text{ if } i= 1) \\
0
& \qquad (\text{ if } i\neq 1) \\
\endcases  \\
R^i \pi_!( K(W;x_0^{\alpha_0}) \otimes \Cal L_{\psi})
& =
\cases
\tilde j_!K(y^{-\alpha_0}) \otimes g(x_0^{\alpha_0}, \psi ) 
& \qquad (\text{ if } i= 1) \\
0
& \qquad (\text{ if } i\neq 1) \\
\endcases  \\
\endalign
$$
and if $\alpha_0$ is trivial, we have
$$
\align
R^i \pi_* \Cal L_{\psi}
& =
\cases
\tilde j_!K(-1) 
& \qquad (\text{ if } i= 1) \\
0
& \qquad (\text{ if } i\neq 1) \\
\endcases  \\
R^i \pi_!\Cal L_{\psi}
& =
\cases
K  
& \qquad (\text{ if } i= 1) \\
\tilde i_* K(-1)  
& \qquad (\text{ if } i= 2) \\
0
& \qquad (\text{ if } i\neq 1,2), \\
\endcases  \\
\endalign
$$
where $\tilde j$ and $\tilde i$ are the natural inclusions
$\bold T(y) \to \bold A(y)$ and $\{ 0\} \to \bold A(y)$.
Note that the variety $\Cal X_D$ is isomorphic to an open
subset $U$ of $\bold T(x) \times \bold T(a)$ defined
by $U = \{ (x, a) \in\bold T(x) \times \bold T(a) \mid
f(a,x) \neq 0\}$. The natural morphism from 
$\bold T(x) \times \bold T(a)$ to $\bold T(a)$ is denoted by 
$\varphi'$.
If $\alpha_0 \neq 0$, then we have
$$
\align
& 
\bold R \varphi^{(p)}_*K(\Cal X^{(p)} ;
x_0^{\alpha_0} \prod_{i=1}^nx_i^{\alpha_i})\otimes \Cal L_{\psi}(t) \\
\simeq &
\bold R \varphi'_*[(
\bold R \pi'_*(F')^*(K(W;x_0^{\alpha_0})\otimes 
\Cal L_{\psi}(t)))\otimes 
K(\bold T(x) \times \bold T(a) ;\prod_{i=1}^nx_i^{\alpha_i})] \\
\simeq &
\bold R \varphi'_*[(F^*\bold R \pi_*(K(W;x_0^{\alpha_0})\otimes
\Cal L_{\psi}(t)))\otimes 
K(\bold T(x) \times \bold T(a) ;\prod_{i=1}^nx_i^{\alpha_i})]  \\
\simeq &
\bold R \varphi'_*(F^*\tilde j_*K(y^{-\alpha_0})\otimes
K(\bold T(x) \times \bold T(a) ;\prod_{i=1}^nx_i^{\alpha_i})) 
\otimes g(x_0^{\alpha_0}, \psi)[-1] \\
\simeq &
\bold R \varphi_*K(\Cal X ;x_0^{\alpha_0}\prod_{i=1}^nx_i^{\alpha_i})
\otimes g(x_0^{\alpha_0}, \psi)[-1]. \\
\endalign
$$
Similarly, we have
$$
\bold R \varphi^{(p)}_!
K(\Cal X^{(p)} ;x_0^{\alpha_0} \prod_{i=1}^nx_i^{\alpha_i})
\otimes \Cal L_{\psi}(t)  \simeq
\bold R \varphi_!K(\Cal X ;x_0^{\alpha_0}\prod_{i=1}^nx_i^{\alpha_i})
\otimes g(x_0^{\alpha_0}, \psi)[-1]. 
$$
If $\alpha_0 = 0$, by using exact sequence and triangle
$$
\align
0 \to j_!j^*K(\bold T(x) \times \bold T(a) ;
 \prod_{i=1}^nx_i^{\alpha_i}) & \to
 K(\bold T(x) \times \bold T(a) ;
 \prod_{i=1}^nx_i^{\alpha_i}) \\
& \to
 i_*i^* K(\bold T(x) \times \bold T(a) ;
 \prod_{i=1}^nx_i^{\alpha_i}) \to 0 \\
\endalign
$$
$$
\align
 K(\bold T(x) \times \bold T(a) ;
 \prod_{i=1}^nx_i^{\alpha_i})[-1] & \to
F^*\bold R\pi_!\Cal L_{\psi} \otimes
K(\bold T(x) \times \bold T(a) ;
 \prod_{i=1}^nx_i^{\alpha_i})  \\
& \to
i_*i^*K(\bold T(x) \times \bold T(a) ;
 \prod_{i=1}^nx_i^{\alpha_i})(-1)[-2] 
\overset{+1}\to\to, \\
\endalign
$$
where $j : U \to \bold T(x) \times \bold T(a)$ and
$i: \bold T(x) \times \bold T(a) -U \to \bold T(x) \times \bold T(a) $
and 
$$
 \bold R \varphi'_* 
K(\bold T(x) \times \bold T(a) ;
 \prod_{i=1}^nx_i^{\alpha_i})(-1) \simeq 0 
\text{ and }
 \bold R \varphi'_!
K(\bold T(x) \times \bold T(a) ;
 \prod_{i=1}^nx_i^{\alpha_i}) \simeq 0, 
$$
we have
$$
\align 
& \bold R\varphi^{(p)}_*K(\Cal X^{(p)} ;\prod_{i=1}^nx_i^{\alpha_i})
\otimes \Cal L_{\psi}(t) \\
\simeq &
\bold R \varphi'_*(\bold R \pi'_* (F')^*\Cal L_{\psi}(t))\otimes
K(\bold T(x) \times \bold T(a) ;
\prod_{i=1}^nx_i^{\alpha_i}) \\
\simeq &
\bold R \varphi'_*(F^*\tilde j_!K)\otimes
K(\bold T(x) \times \bold T(a) ;
\prod_{i=1}^nx_i^{\alpha_i})(-1)[-1] \\
\simeq &
\bold R \varphi'_*(j_!j^*
K(\bold T(x) \times \bold T(a) ;\prod_{i=1}^nx_i^{\alpha_i}))(-1)[-1] \\
\simeq &
\bold R \varphi'_*(i_*i^*
K(\bold T(x) \times \bold T(a) ;\prod_{i=1}^nx_i^{\alpha_i}))(-1)[-2] \\
\endalign
$$
and
$$
\align 
& \bold R\varphi^{(p)}_!K(\Cal X^{(p)} ;\prod_{i=1}^nx_i^{\alpha_i})
\otimes\Cal L_{\psi}(t ) \\
\simeq &
\bold R \varphi'_!(\bold R \pi'_! (F')^*\Cal L_{\psi}(t))\otimes
K(\bold T(x) \times \bold T(a) ;\prod_{i=1}^nx_i^{\alpha_i}) \\
\simeq &
\bold R \varphi'_!F^*(\bold R \pi_!\Cal L_{\psi}(t))\otimes
K(\bold T(x) \times \bold T(a) ;\prod_{i=1}^nx_i^{\alpha_i}) \\
\simeq &
\bold R \varphi'_!(i_*i^*
K(\bold T(x) \times \bold T(a) ;\prod_{i=1}^nx_i^{\alpha_i}))(-1)[-2]. \\
\endalign
$$
Therefore we have the proposition.
\enddemo

As a consequence, we have the following corollary.
\proclaim{Corollary}
If $k=\bold C$ (resp. $k=\bold F_q$),
a hypergeometric sheaf $\Cal G(D)$ (resp. $\Cal G(D, \psi )$) 
is an irreducible perverse sheaf.
\endproclaim
\demo{Definition}
Let $\bold T$ be a torus over $k=\bold C$ (resp. $k=\bold F_q$).  
A local system $\Cal F$ 
on an open set $U$ in $\bold T$
is called hypergeometric sheaf 
if there is an isomorphism 
$\phi: \bold T \to \bold T(R)$ and a non-resonant hypergeometric data $D$
such that the restriction $\phi^*\Cal G(D)\mid_U$ 
(resp. $\phi^*\Cal G(D, \psi )\mid_U$) of
$\phi^*\Cal G(D)$ (resp. $\phi^*\Cal G(D, \psi )$)
to $U$ is isomorphic to $\Cal F$ on $\bold T$.
\enddemo

\heading
\S 3 Variation of $K$-Hodge structure
\endheading
Throughout this section, we assume that $k=\bold C$.
First we give the definition of variation of Hodge structures.
Let $S$ be a complex manifold and $n$ be a non-negative integer.
A variation of $K$-Hodge structure is a datum 
$(V,(F^{\bullet}_{\sigma}))$ consisting of
\roster
\item
a local system $V$ of $K$-vector spaces on $S$ and
\item
a decreasing filtration $F^{\bullet}_{\sigma}$ by locally
free $\Cal O_S$-module on $V\otimes_{K,\sigma}\Cal O_S$
for all the embeddings $\sigma : K \to \bold C$.
\endroster
For this data, we insist the following conditions on 
$F^{\bullet}_{\sigma}$;
\roster
\item
For any point $s \in S$, the natural homomorphism
$$
(F^i_{\sigma})_s \oplus (1\otimes c)(F^{n-i}_{\bar\sigma})_s \to
V_s \otimes_{K,\sigma}\bold C
$$
is an isomorphism, where $\bar\sigma$ is the complex conjugate of
the embedding $\sigma$,
$c$ the complex conjugate, and $1\otimes c$
the natural homomorphism 
$$
1 \otimes c :V_s 
\otimes_{K,\bar\sigma}\bold C \to V_s \otimes_{K,\sigma}\bold C.
$$
\item
(Griffith transversality) Let $\nabla_{\sigma}$ be the connection on
$V\otimes_{K,\sigma}\Cal O_S$ defined by 
$1\otimes d : 
V\otimes_{K,\sigma}\Cal O_S \to V\otimes_{K,\sigma}\Omega^1_S.$
Then $\nabla_{\sigma}(F^i_{\sigma}) \subset 
F^{i-1}_{\sigma}\otimes \Omega^1_S$.  
\endroster

Let $D=(R, \{l_i\}_i, \{\kappa_i\}_i)$ be a hypergeometric data
and $U_R$ be the open set of $\bold T(R)$ consisting of non-degenerate
points.  We introduce a variation of $K$-Hodge structure on the
local system $\Cal G(D)\mid_{U_R}$.
For simplicity, we assume $\sum_{i=1}^r\kappa_i \neq 0$. For the case 
$\sum_{i=1}^r\kappa_i =0$, we can define the variation of Hodge structure
in the same way.
Let $m$ be the minimal common denominator of $\kappa$. Then we have
$\kappa \in \bold Z^r\otimes I_m(\bold C, K)$.
The super lattice $\frac{1}{m} L$ of $ L$ is denoted by
$ M$. We identify the polytope $\Delta$
in $L\otimes \bold R$ with that in $M\otimes \bold R$.
For a non-negative integer $i$, the $\bold C$ subspace of $\bold C[ M]$ 
generated by the set of monomials
$\{ m \in M \mid  m \in i\Delta \}$ is denoted by
$L_M(i\Delta)$. Let $u$ be a indeterminacy. The ring 
$R = \oplus_{i \geq o}L_M(i\Delta )u^i$ is a graded ring 
with $\deg (u)=1$. The projective variety $\Proj (R)$ is denoted
by $\bold P_{\Delta}$. This is a toric variety associated to 
the fan $F_{\Delta}=\{C_b\}_{b\in \Delta}$, where
$$
C_b=\{ x\in (L)^*\otimes \bold R \mid
(x, b'-b) \geq 0 \text{ for all } b' \in \Delta \}
$$
We take a regular triangulation $\tilde F$ of the fan
$F_{\Delta}$ with respect to the lattice $M$.
Then the associate toric variety
$X(\tilde F)$ is smooth. The closure of the pull back of
the divisor $\Cal X \in \bold T( L) \times \bold T(a)$
in $X(\tilde F)$ via the morphism 
$\bold T( M)\times \bold T(a) \to 
\bold T( L)\times \bold T(a)$ is denoted by $\bar\Cal X_M$.
By the definition of non degeneracy, the restriction 
$\varphi_M : \bar\Cal X_M \mid_{U_a} \to U_a$ of the morphism
$\bar\Cal X_M \to \bold T(a)$ to the open set $U_a$ is a smooth.
Therefore $R^i{\varphi_M}_* K$ is a variation of $K$-Hodge structure on $U_a$.
By the construction, the covering group 
$\Aut(\bold T(M)/\bold T(L))\simeq 
\Hom((M/L),\mu_m(\bold C))$ acts on the variety $\bar\Cal X_M$.
The element $\alpha \in L\otimes I(\bold C, K)$
defines a character of $\Aut (\bold T (M)/\bold T(L))$.
By the non-resonance condition, 
$R^i\varphi_*K(\Cal X;x^\alpha )\mid_{U_a}$ is naturally 
isomorphic to the $\alpha$ part 
$(R^i{\varphi_M}_*K)(\alpha )$ of $R^i{\varphi_M}_*K$.  
It defines a variation
of $K$-Hodge structure on $R^i\varphi_*K(\Cal X;x^\alpha )\mid_{U_a}$.

\heading
\S 4 Multiplicative equivalence
\endheading

In this section, we fix a lattice $R$.
Let $I = [1, r]$ and $D = (R, \{l_i\}_{i\in I},\{\kappa_i\}_{i \in I})$ 
be a hypergeometric data. A permutation $\sigma$ of the index $I$,
gives another hypergeometric data 
$D_{\sigma}=(R, \{l_{\sigma (i)}\}_i, \{\kappa_{\sigma (i)}\}_i)$.
The sheaf $\Cal G(D)$ and $\Cal G(D_{\sigma})$ are isomorphic
to each other.
Therefore we use the notation $D=(l_1,\kappa_1) \oplus \cdots \oplus
(l_r, \kappa_r)=\oplus_{i=1}^r(l_i.\kappa_i)$.
Notice that $(l_i,\kappa_i)$ is not necessarily a hypergeometric 
data.
More generally, for data $(\{l_i\}_{i \in I},\{\kappa_i\}_{i\in I})$
and $(\{l'_j\}_{j \in J},\{\kappa'_j\}_{j\in J})$,
$(\{l_i\}_{i \in I}\cup\{l'_j\}_{j \in J},
\{\kappa_i\}_{i\in I}\cup\{\kappa'_j\}_{j\in J})$
is expressed as 
$(\{l_i\}_{i \in I},\{\kappa_i\}_{i\in I}) \oplus
(\{l'_j\}_{j \in J},\{\kappa'_j\}_{j\in J})$ and so on.
We use the above notation even though each terms do not satisfy
the conditions for hypergeometric data.
\demo{Remark}
The notation $\oplus$ corresponds to convolution product in [G-L].
\enddemo
If $l_1=\cdots = l_r=l$, 
$\{l_i\}_{i=1, \dots r}$ is denoted by $l^r$.
Let $d$ be a positive natural number prime to the 
characteristic $p$.
For an element $\kappa \in I(\bar k, K)$, 
the inverse image of the multiplication by $d$ on $I(\bar k,K)$  
is denoted by $\frac{1}{d}\{\kappa\}$. The cardinality of this set is 
$d$.  Let $D = (R, \{l_i\}_i.\{\kappa\}_i)$ be a hypergeometric data.
Suppose that $l_i$ is not equal to 0 and divisible by $d$.
We define the hypergeometric data $D^{(i,d)}$ by
$$
D^{(i,d)} = \oplus_{j \neq i}(l_j,\kappa_j) \oplus 
((\frac{l_i}{d})^d, \frac{1}{d}\{\kappa_i\})
$$

By the following lemma, $D^{(i,d)}$ satisfies the conditions of
hypergeometric data.
\proclaim{Proposition 4.1}
The data $D^{(i,d)}$ is separated.  If $D$ is non-resonant, then 
$D^{(i,d)}$ is also non-resonant.
\endproclaim
\demo{Proof}
We assume $d >1$, $i=1$. 
If $(l_1/d)(r)=1$, this $r$ does not give non separatedness criterion.
If $l_i(r)=1$ and 
$l_j(r)=1$ and $l_k(r)=0$ for $r \neq 1,i,j$ and $l_1(r)/d =0$.
Then this contradicts to the separatedness of the data $D$.
Before proving the resonance of $D^{(1,d)}$, we give a description of
$L^{(1,d)} = 
\Coker(\bold l^{(1,d)}:R \to \bold Z^d \oplus \bold Z^{r-1})$,
where $\bold l^{(1,d)}=(l^{(1,d)}_i)_i$ is given by the data 
$D^{(1,d)}=(R, \{l^{(1,d)}_i\}_i,\{\kappa^{(1,d)}_i\}_i)$.
Let $q^{(1,d)}$ be the natural projection $\bold Z^d \oplus \bold Z^{r-1}$
and the image of $((0, \dots ,\overset{i}\to{1} ,\dots ,0),(0,\dots ,0))$
($i=1, \dots ,d$)
and $((0,\dots ,0),(0, \dots ,\overset{j}\to{1} ,\dots ,0))$
($j=2, \dots ,r$)
under the homomorphism $q^{(1,d)}$ is denoted by $\eta_i$ and
$\omega_j$ respectively.
By the definition, $L^{(1,d)}$ is generated by 
$\eta_i$ and $\omega_j$
($i=1, \dots ,d, j=2, \dots ,r$).
By the following commutative diagram,
the submodule $Q$ of $L$ generated by 
$\omega_i$ ($i=2, \dots ,r$) 
and $d\omega_1$ is isomorphic to the submodule
of $L^{(1,d)}$ generated by $\omega_i$ ($i=2, \dots ,r$)
and $\eta_1+ \cdots +\eta_d$.
$$
\CD
& R \\
 \bold l^{(1,d)} \swarrow & & \searrow \bold l \\
\bold Z^d \oplus \bold Z^{r-1} &  
@VV{(\frac{l_1}{d},l_2,\dots,l_r)}V \bold Z 
\oplus \bold Z^{r-1} \\
 (\Delta ,1) \nwarrow & & \nearrow  (d,1) \\
& \bold Z \oplus \bold Z^{r-1} \\
\endCD
$$
Therefore $L^{(1,d)}$ is isomorphic to
$(Q \oplus \oplus_{i=1}^d \bold Z \eta_i)/
(\eta_1+ \dots + \eta_d - d\omega_1)\bold Z$. 
Then by the definition of
$\alpha = q(\kappa )$ and $\alpha^{(1,d)}=q^{(1,d)}(\kappa^{(1,d)})$, 
$$
\align
\alpha = & \kappa_1\omega_1 + \sum_{i=2}^r \kappa_i \omega_i \\
\alpha^{(1,d)} = 
& \sum_{i=1}^d\frac{\kappa_1+i}{d}\eta_i + 
\sum_{i=2}^r \kappa_i \omega_i. \\
\endalign
$$
Let $\Delta$ and $\Delta^{(1,d)}$ be a convex hull of  
$\{\omega_i\}_{i=1. \dots , r}\cup \{0\}$ and
$\{\eta_i\}_{i=1, \dots ,d}\cup \{\omega_j\}_{j=2, \dots ,r}
\cup \{0\}$ in $L \otimes \bold R$ and 
$L^{(1,d)} \otimes \bold R$ respectively. Since the barry center
$(\eta_1 + \cdots + \eta_d)/d = \omega_1$ 
of $\eta_1, \dots ,\eta_d$ is contained in 
$\Delta$, $\Delta^{(1,d)}$ is a convex hull of $\Delta$
and $\eta_1, \dots ,\eta_d$.

Therefore the linear hull of a codimension 1 
face of $C(\Delta^{(1,d)})$ is the linear hull $H_{i,\sigma}$ of
$\eta_j$ ($j \neq i$)
and a codimension 1 face $\sigma$ of $C(\Delta)$.
Let $h_{\sigma} $ be a primitive linear form
vanishing on $\sigma$.
Let $\eta_i^*$ be the linear form on 
$Q \oplus \oplus_{i=1}^d \bold Z \eta_i$
such that $\eta_i^*(Q)=0$ and 
$\eta_i^*(\eta_j)=\delta_{i,j}$.
The linear form defined by
$$
h_{i,\sigma}=dh_{\sigma}(\omega_1)\eta_i^* +h_{\sigma}
$$
vanishes on $\sigma$ and $\eta_j$ ($j \neq i$). 
Since
$h_{i,\sigma}(d\omega_1-\sum_{i=1}^d\eta_i)=0$,
it defines a linear form on $L^{(1,d)}$
vanishing on the hyperplane $H_{i,\sigma}$. 
Since 
$$
h_{i,\sigma}(\alpha^{(1,d)} ) \equiv 
h_{\sigma}(\alpha ) \quad
(\text{mod } \bold Z),
$$
and the linear form
$h_{i,\sigma}$ is an integral on $L^{(1,d)}$, 
$\alpha^{(1,d)}$ is non-resonant.
\enddemo
\demo{Definition}
The equivalence relation generated by $D \sim D^{(1,d)}$ is called
multiplicative equivalence.  
\enddemo

\heading
\S 5 Constant equivalence and algebraic correspondences
\endheading
\heading
\S 5.1 Reduced part and constant correspondence for $k=\bold C$
\endheading

Let $D = (R, \{l_i\}_i,\{\kappa_i\}_i)$ be a hypergeometric data and
$S= \{ i \mid l_i =0\}$.  Then it is easy to see that the homomorphism
$R \to \oplus_{i \notin S}\bold Z$ is also a primitive separate
embedding.  Let $L_{red} = \Coker (R \to \oplus_{i \notin S}\bold Z)$.
Then we have $L = \Coker (R \to \bold Z^r) \simeq
L_{red} \oplus_{i \in S} \bold Z \omega_i$.
\proclaim{Proposition 5.1}
If $D$ is non-resonant, then $D_{red} = 
(R, \{l_i\}_{i \notin S}, \{\kappa_i\}_{i \notin S})$ is non-resonant
and $\kappa_i \notin \bold Z$ ($i\in S$).
\endproclaim
\demo{Proof}
Let $C(\Delta_{red})$ be the convex cone generated by
of $\omega_i$ ($ i \notin S$). Then the codimension 1 face
of the convex cone $C(\Delta )$ generated by
$\Delta$ is either
\roster
\item
the cone generated by $C(\Delta_{red})$ and $\omega_i$ ($i \in S -\{k\}$) 
for some $k \in S$ or
\item
the cone generated by $\sigma$ and
$\omega_i$ ($i \in S$), where $\sigma$ is a codimension
1 face of $C(\Delta_{red})$.
\endroster 
The defining equation $\tilde h_{\sigma}$ of the linear hull of
$\sigma$ and $\omega_i$,($i \in S$) is the same as the equation
$h_{\sigma}$ of that of $\sigma$ under the decomposition.
Therefore 
$h_{\sigma}(\alpha_{red} ) = \tilde h_{\sigma}(\alpha )$, where
$\alpha_{red} = \sum_{i\notin S}\kappa_i \omega_i$.
The defining equation of the hyperplane of type (1) is 
the dual base
$\omega^*_k$ ($k \in S$) of $\{\omega_i\}_{i\in S}$ and $L_{red}$.
Therefore $\kappa_k \notin \bold Z$.
\enddemo
\demo{Definition}
The above hypergeometric data $D_{red}$ is called the reduced part
of $D$.
On the set of hypergeometric data over a lattice $R$, we introduce
the constant equivalence as
the equivalence relation generated by $D \sim D_{red}$.
\enddemo
We define Fermat motif and Artin-Shreier motif as follows.
Let $\kappa_1, \dots , \kappa_r$ be an element of 
$I(k, K)-\{0\}$.  Suppose that $\sum_{i=1}^r\kappa_i \notin \bold Z$
if $k =\bold C$.  Positive Fermat motif (resp. 
Positive Fermat-Artin-Shreier motif) is 
a $K$-Hodge structure ($l$-adic representation ) of the form
$H^{r-1}(F, K(\kappa_1, \dots ,\kappa_r))$
(resp. $H^n(FAS,K(\kappa_1, \dots , \kappa_r))(\psi)$),
where $F$ and $FAS$ is defined by
$$
\align
F: x_1 + \cdots + x_r = & 1 \\
FAS: x_1 + \cdots + x_r = &\tau^p -\tau \\
\endalign
$$
and the local system $K(\kappa_1, \dots ,\kappa_r)$ 
on $F$ and $FAS$ correspond to the Kummer character
$\kappa (\prod_{i=1}^r x_i^{\kappa_i})$.
We omit the expression for dimension of $F$ and $FAS$,
because it is determined by the 
number of the components of the Kummer characters.
It is well known that they are stable under the tensor product;
$$
\align
& H^{r-1}(F)(\kappa_1, \dots, \kappa_r) \otimes 
H^{s-1}(F)(\lambda_1, \dots ,\lambda_s)(-1) \\
\simeq & 
H^{r+s}(F)(\kappa_1, \dots, \kappa_r,\lambda_1, \dots ,\lambda_s,
-\sum_{i=1}^r\kappa_i) \otimes K((-1)^{-\sum_{i=1}^r\kappa_i})
 \\
& H^{r}(FAS)(\kappa_1, \dots, \kappa_r, \psi) \otimes 
H^{s}(FAS)(\lambda_1, \dots ,\lambda_s, \psi) \\
\simeq & 
H^{r+s}(FAS)(\kappa_1, \dots, \kappa_r,\lambda_1, \dots ,\lambda_s,
\psi),  \\
\endalign
$$
where $K((-1)^{-\sum_{i=1}^r\kappa_i})$ is the Hodge structure on
$\Spec \bold C$ corresponding to the Kummer character
$(-1)^{-\sum_{i=1}^r\kappa_i}$.  As Hodge structures,
$K((-1)^{-\sum_{i=1}^r\kappa_i})$ is isomorphic to $K$ over $\Spec \bold C$.
Moreover these isomorphisms are induced by algebraic correspondences.
Positive Fermat motif contains 
$K(-1) \simeq H(F)(\kappa_1,-\kappa_1,\kappa_2)$ and Fermat motif is
generated by tensoring positive Fermat motif and $K(1)$.  Using
$K(-1) \simeq FAS(-\kappa_1, \kappa_1)$, we define Fermat-Artin-Shreier
motif in the same way.
An element in $H(F)(\kappa_1, \dots ,\kappa_r)$  is called Hodge cycle
if the Hodge type of
$H(F)(\kappa_1, \dots ,\kappa_r)$ is $(m,m)$ for all the embedding
$\sigma : K \subset \bold C$.
An element in
$H(FAS)(\kappa_1, \dots ,\kappa_r,\psi)$ 
is called Tate cycle if the action of Frobenius of $\bold F_q$ is $q^m\times \zeta$,
where $\zeta$ is a root of unity.  In many case, Hodge cycles and
Tate cycles of Fermat hypersurface is know to be algebraic cycle.
But it remains still open in general case ([A][S]).

We introduce an algebraic correspondence associated to
$\Cal G(D)$ (resp. $\Cal G(D, \psi )$) and 
$\Cal G(D_{red})$ (resp $\Cal G(D_{red}, \psi )$).
First we consider the case $k=\bold C$.
\proclaim{Theorem 1} Let $D$ and $D_{red}$ be hypergeometric
data defined in Proposition 5.1. We renumber the index such that
$S = \{ 1, \dots ,s\}$.  Then there exists an isomorphism of variation of 
Hodge structures:
$$
\align
\Cal G(D)  \simeq &\Cal G(D_{red}) \otimes H^s(F)
(\kappa_1 ,\dots , \kappa_s,\sum_{i=s+1}^r\kappa_i ) 
\qquad (\sum_{i=1}^r\kappa_i \neq 0, \sum_{i=s+1}^r\kappa_i \neq 0) \\
\Cal G(D) \simeq & \Cal G(D_{red}) \otimes H^s(F)
(\kappa_1 ,\dots , \kappa_s)(-1) 
\qquad (\sum_{i=1}^r\kappa_i \neq 0, \sum_{i=s+1}^r\kappa_i = 0) \\
\Cal G(D) \simeq & \Cal G(D_{red}) \otimes H^s(F)
(\kappa_1 ,\dots , \kappa_s) \otimes
K((-1)^{\sum_{i=1}^s\kappa_i})  \\
& \qquad (\sum_{i=1}^r\kappa_i = 0, \sum_{i=s+1}^r\kappa_i \neq 0) \\
\Cal G(D) \simeq & \Cal G(D_{red}) \otimes H^s(F)
(\kappa_1 ,\dots , \kappa_{s-1})(-1)\otimes
K((-1)^{\sum_{i=1}^{s-1}\kappa_i}) \\
& \qquad (\sum_{i=1}^r\kappa_i = 0, \sum_{i=s+1}^r\kappa_i = 0) \\
\endalign 
$$
defined by an algebraic correspondence.
\endproclaim
We use systems of 
coordinates $y = (y_1, \dots , y_s)$ and 
$\bar x=(x_{s+1},\dots ,x_{n})$.
We take a base of $L$ such that 
$x^{\omega_i}=x_0y_i$ for $i=1, \dots ,s$ and 
$x^{\omega_i}=x_0\bar x^{\bar\omega_i}$ ($i=s+1, \dots , n$).  
For a constructible sheaf $\Cal F$ on a torus $\bold T$ 
is said to be of weight less (resp. more ) than $w$ if
there exists an open set $U$ of $\bold T$ such that the fiber of
$\Cal F$ at $p$ in $U$ is mixed Hodge structure of weight 
less (resp. more ) than $w$ and it is denoted by
$\wt (\Cal F) \leq w$ (resp. $\wt (\Cal F) \geq w$).
\heading 
Case (a) $\sum_{i=1}^r\kappa_i \neq 0$
\endheading
The variety $\Cal X_D$ associated to $D$ and the 
rank 1 local system is given by
$$
\align
\Cal X_D = \{ &  (y,x_0,\bar x,a) \in \bold T(y) \times \bold T(x_0) 
\times \bold T(\bar x) \times \bold T(a) \\
& \mid
x_0(a_1 y_1 + \cdots a_s y_s + \sum_{i=s+1}^ra_i\bar x^{\bar\omega_i})=1 \} \\
\endalign
$$ 
and 
$$
K(x^\alpha a^\kappa )=
K (\Cal X_D ; \prod_{i=1}^s (x_0a_iy_i)^{\kappa_i}
\prod_{i\geq s+1}(x_0a_i\bar x^{\bar\omega_i})^{\kappa_i})
$$
We put $f_{0}=\sum_{i=s+1}^ra_i \bar x^{\bar\omega_i}$.
Let $\Cal X_D^0$ and $\Cal X_D^1$ be the open set and
closed set of $\Cal X_D$ defined by
$$
\align
\Cal X_D^0 = &\{ (y,x_0, \bar x,a) \in \Cal X_D \mid 
f_0 \neq 0 \} \text{ and } \\
\Cal X_D^1 = &\{ (y,x_0,\bar x ,a) \in \Cal X_D \mid 
f_0 = 0 \}  \\
\endalign
$$
Let $\varphi^0$ (resp. $\varphi^1$) be the composite of
the open immersion $\Cal X^0_D \to \Cal X_D$ (resp.
the closed immersion $\Cal X^1_D \to \Cal X_D$) and 
$\varphi :\Cal X_D \to \bold T(a)$. 
We get the following long exact sequence
$$
\align 
\cdots & \to R^n{\varphi_{\Cal X^1}}_*K(x^\alpha a^\kappa) 
\overset{i_*i^{!}}\to\to
R^n\varphi_* K(x^\alpha a^\kappa) \to
R^n{\varphi^0}_* K(x^\alpha a^\kappa)   \\
& \overset{\delta}\to\to  
R^{n+1}{\varphi_{\Cal X^1}}_*K(x^\alpha a^\kappa)
\simeq R^{n-1}\varphi^1_*K(x^\alpha a^\kappa)(-1) \to \cdots, \\
\endalign
$$
where $R^i{\varphi_{\Cal X^1}}_*$ denotes the $i$-th higher direct image
with a support in $\Cal X^1_D$.  
\heading
Case (a-1) $\sum_{i=1}^r\kappa_i \neq 0$ and $\sum_{i=s+1}^r\kappa_i \neq 0$
\endheading
Let us introduce sets of variables $Y =(Y_1, \dots ,Y_s)$, $X_0$,
and $t_0$ and
change coordinate by $Y_i = x_0a_iy_i$, $X_0 = x_0f_0$ and
$t_0 f_0 =1$. Then $\Cal X_D^0$ is isomorphic to $\Cal X_{D_{red}} \times F$,
where
$$
\align
\Cal X_{D_{red}} = &\{ (\bar x,t_0, a) \in  
\bold T(\bar x)\times \bold T(t_0) \times \bold T(a) \mid
t_0f_0(a,\bar x)=1 \} \\
F= &\{ (Y,X_0) \in \bold T(Y) \times \bold T(X_0) \mid 
Y_1 + \dots + Y_s + X_0 =1 \} \\
\endalign
$$
Under the above isomorphism, the corresponding rank 1 local system is
$pr_1^*(\Cal F^{(1)}) \otimes pr_2^*(\Cal F^{(2)})$,
where
$$
\Cal F^{(1)} =
K (\Cal X_{D_{red}}; 
\prod_{i=s+1}^r(a_i\bar x^{\bar\omega_i})^{\kappa_i}
t_0^{\sum_{i=s+1}^r\kappa_i})
\text{ and }
\Cal F^{(2)} =
K (F; \prod_{i=1}^sY_i^{\kappa_i}X_0^{\sum_{i=s+1}^r\kappa_i}).
$$
The natural homomorphism $\Cal X_{D_{red}} \to \bold T(a_{s+1}, \dots , a_r)$
is denoted by $\varphi_{red}$. 
Therefore we have 
$$
R^i \varphi^0_* K(x^\alpha a^\kappa)
\simeq 
R^{n-s}\varphi_{red *}\Cal F^{(1)} 
\otimes H^s(F)(\kappa_1, \dots ,\kappa_s, 
\sum_{i=s+1}^r \kappa_i).
$$ 
We have
$\wt (R^{n}{\varphi^0}_*K(x^\alpha a^\kappa))=n$ and
$\wt (R^{n+1}{\varphi_{\Cal X^1}}_*K(x^\alpha a^\kappa)) \geq n+1$.
Therefore the connecting homomorphism $\delta$
vanishes. Since both $R^n\varphi_* K(x^\alpha a^\kappa)$ and
$R^n{\varphi^0}_* K(x^\alpha a^\kappa)$ are perverse and irreducible,
we get the required isomorphism.
\heading
Case (a-2) $\sum_{i=1}^r\kappa_i \neq 0$ and $\sum_{i=s+1}^r\kappa_i = 0$
\endheading
We introduce an coordinate $Y=(Y_1, \dots , Y_s)$ and
change coordinates 
\linebreak
$(Y,x_0,\bar x,a) \to (y, x_0, \bar x,a)$ by
$Y_i = x_0a_iy_i$.
Then the variety $\Cal X_D^1$ is isomorphic to 
$\Cal X_{D_{red}} \times F \times \bold T( x_0)$, where
$$
\align
\Cal X_{D_{red}}= & \{\bar x \in \bold T(\bar x) 
\mid f_0 = 0 \}\\
F= &  \{ Y \in \bold T(Y) \mid Y_1 + \cdots Y_s =1 \}. \\
\endalign
$$
The corresponding rank 1 local system is  
$pr_1^* (\Cal F^{(1)}) \otimes pr_2^* (\Cal F^{(2)})$,
where $\Cal F^{(1)}$ and $\Cal F^{(2)}$ is the local system
on $\Cal X_{D_{red}}$ and $F$
defined by
$$
\Cal F^{(1)} =
K(\Cal X_{D_{red}}; \prod_{i=s+1}^r (a_i\bar x^{\bar\omega_i})^{\kappa_i})
\text{ and }
\Cal F^{(2)} =
K(F ; \prod_{i=1}^sY_i^{\kappa_i}).
$$  Therefore we have
$$
R^{n}\varphi_{\Cal X^1 *}K(x^\alpha a^\kappa) =
R^{n-s-1}{\varphi_{{red *}}}\Cal F^{(1)} \otimes
H^{s-1}(F)(\kappa_1, \dots ,\kappa_s)(-1).
$$
Here $\varphi_{red}$ denotes the natural morphism 
$\Cal X_{D_{red}} \to \bold T(a_{s+1},\dots ,a_r)$.

We show that $\wt (R^n \varphi^0_* K(x^\alpha a^\kappa )) =n+2$. In 
fact, under the change of coordinate $Y_i = x_0 a_i y_i$
and $X_0 = x_0 f_0$, the variety $\Cal X^0_D$ is isomorphic to 
$U_1 \times U_2 \times \bold T(x_0)$, where
$$
\align
U_1 = &\{ Y \in \bold T(Y) \mid Y_1+ \dots Y_s - 1 \neq 0\} \\
U_2 = &\{ \bar x \in \bold T(\bar x) \times  \bold T(a) 
\mid f_0(a, \bar x) \neq 0\} \\
\endalign
$$
and the corresponding Kummer character is $pr_1^*\Cal F^{(2)} \otimes
pr_2^*\Cal F^{(1)}$ defined as above.  Therefore
$\wt (H^s(U_1, \Cal F^{(2)}))= s+1$ and 
$\wt (R^{n-s}{\pi_2}_* \Cal F^{(1)})= n-s+1$, where 
$\pi_2 : U_2 \to \bold T(a_{s+1},\dots ,a_r)$ 
is the natural projection. Therefore the 
morphism $i_*i^{!}$ is surjective on the open set $U_a$.
The perversity of $R^n\varphi_{\Cal X^1 *} K(x^\alpha a^\kappa )$
and $R^n\varphi_* K(x^\alpha a^\kappa )$ implies the required isomorphism.
\heading
Case (b) $\sum_{i=1}^r\kappa_i = 0$ 
\endheading
In this case, the variety $\Cal X_D$ 
and the rank 1 local system
associated to 
$D$ are given by
$$
\Cal X_D = \{ (y,\bar x,a) \in \bold T(y) \times  
 \bold T(\bar x) \times \bold T(a)
\mid
a_1 y_1 + \cdots a_s y_s + \sum_{i=s+1}^ra_i\bar x^{\bar\omega_i}=0 \}
$$ 
and 
$$
K(x^\alpha a^\kappa)=
K (\Cal X_D ; \prod_{i=1}^s (a_iy_i)^{\kappa_i}
\prod_{i\geq s+1}(a_i\bar x^{\bar\omega_i})^{\kappa_i}).
$$
We put $f_0=\sum_{i=s+1}^ra_i \bar x^{\bar\omega_i}$.
Let $\Cal X_D^0$ and $\Cal X_D^1$ be the open set and
closed set of $\Cal X_D$ defined by
$$
\Cal X_D^0 = \{ (y, \bar x,a) \in \Cal X_D \mid 
f_0 \neq 0 \} \text{ and } 
\Cal X_D^1 = \{ (y,\bar x ,a) \in \Cal X_D \mid 
f_0 = 0 \}  
$$
Let $\varphi^0$ (resp. $\varphi^1$) be the composite of
the open immersion $\Cal X^0_D \to \Cal X_D$ (resp.
the closed immersion $\Cal X^1_D \to \Cal X_D$) and $\varphi$.
On the other hand, we get the following long exact sequence
$$
\align
\cdots & \to R^{n-1}{\varphi_{\Cal X^1}}_*K(x^\alpha a^\kappa) 
\overset{i_*i^{!}}\to\to
R^{n-1}\varphi_* K(x^\alpha a^\kappa) \to
R^{n-1}{\varphi^0}_* K(x^\alpha a^\kappa)   \\
& \overset{\delta}\to\to  
R^{n}{\varphi_{\Cal X^1}}_*K(x^\alpha a^\kappa)
\simeq R^{n-2}\varphi^1_*K(x^\alpha a^\kappa)(-1) \to \cdots, \\
\endalign
$$
where $R^i{\varphi_{\Cal X^1}}_*$ denotes the $i$-th higher direct image
with a support in $\Cal X^1$.  
\heading
Case (b-1) $\sum_{i=1}^r\kappa_i = 0$ and $\sum_{i=s+1}^r\kappa_i \neq 0$ 
\endheading
We introduce systems of coordinates $t_0=(t_0)$ and $Y = (Y_1, \dots ,Y_s)$.
By changing coordinate, $t_0 f_0 =1$ $Y_i = -t_0 a_i y_i$, the variety
$\Cal X^0_D$ is isomorphic to the product $\Cal X_{D_{red}} \times F$,
where
$$
\align 
\Cal X_{D_{red}} = & \{ (t_0, \bar x,a) \in \bold T(t_0) \times
\bold T (\bar x) \times \bold T(a) \mid t_0 f_0 =1 \} \\
F = & \{ y \in \bold T(Y) \mid \sum_{i=1}^s Y_s = 1 \}, \\
\endalign
$$
and under isomorphism, the corresponding rank 1 local system is
$pr_1^*\Cal F^{(1)} \otimes pr_2^*\Cal F^{(2)}$, where
$$
\Cal F^{(1)} = K (\Cal X_{D_{red}}; \prod_{i=s+1}^r 
(t_0a_i\bar x^{\bar\omega_i})^{\kappa_i}) \text{ and }
\Cal F^{(2)} = K (F; \prod_{i=1}^s (-Y_i)^{\kappa_i}).
$$
The natural morphism $\Cal X_{D_{red}} \to \bold T(a_{s+1},\dots ,a_r)$
is denoted by $\varphi_{red}$.
Therefore we have
$$
R^{n-1}\varphi^0_*K(x^\alpha a^\kappa) =
R^{n-s}{\varphi_{D{red}}}_*\Cal F^{(1)} \otimes
H^{s-1}(F)(\kappa_1, \dots ,\kappa_s)\otimes
K((-1)^{\sum_{i=1}^s\kappa_i}).
$$
Since $\wt (R^n{\varphi_{\Cal X^1}}_* K(x^\alpha a^\kappa )) \geq n$,
the morphism $\delta$ is zero on the open set $U_a$.
Therefore $R^{n-1}\varphi_*K(x^\alpha a^\kappa ) \to
R^{n-1}\varphi^0_* K(x^\alpha a^\kappa )$ is an isomorphism by the 
perversity and irreducibility argument.
\heading
Case (b-2) $\sum_{i=1}^r\kappa_i = 0$ and $\sum_{i=s+1}^r\kappa_i = 0$
\endheading
We introduce a variable $\eta =(\eta_1, \dots ,\eta_{s-1})$
and $Y_s = (Y_s)$.
By changing coordinate, $(y,\bar x,a) \to (\bar x,a, \eta , Y_s)$
with $\eta_i =-\frac{a_iy_i}{a_sy_s}$ ($i=1, \dots , s-1$),
the variety $\Cal X^1_D$ is 
isomorphic to $\Cal X_{D_{red}} \times F \times \bold T(Y_s)$, where
$$
\align
\Cal X_{D_{red}}= & 
\{(\bar x,a)\in \bold T(\bar x) \times \bold T(a) \mid f(a,\bar x)=0 \} \\
F = & 
\{ \eta \in \bold T(\eta) \mid \eta_1 +\cdots + \eta_{s-1} =1 \}, \\
\endalign
$$
and the rank 1 local system is
$pr_1^*\Cal F^{(1)} \otimes pr_2^*\Cal F^{(2)}$, where
$$
\Cal F^{(1)} = K (\Cal X_{D_{red}}; \prod_{i=s+1}^r 
(a_i\bar x^{\bar\omega_i})^{\kappa_i}) \text{ and }
\Cal F^{(2)} = K (F; \prod_{i=1}^{s-1} (-\eta_i)^{\kappa_i}).
$$
The natural map $\Cal X_{D_{red}} \to \bold T(a_{s+1},\dots ,a_r)$
is denoted by $\varphi_{red}$.
Therefore we have
$$
\align
R^{n-1}\varphi_{\Cal X_1,*}K(x^\alpha a^\kappa) =
R^{n-s-1}{\varphi_{D{red},*}}\Cal F^{(1)} & \otimes
H^{s-2}(F)(\kappa_1, \dots ,\kappa_{s-1})(-1) \\
& \otimes
K((-1)^{\sum_{i=1}^{s-1}\kappa_i}). \\
\endalign
$$
As in the case (a-2), we can show that 
$\wt (R^{n-1}\varphi^0_*K(x^{\alpha}a^{\kappa})) \geq n+1$.
Therefore $i_*i^{!}$ is surjective. Perversity and irreducibility
implies the required isomorphism.
We complete all the cases.

\heading
\S 5.2 Constant correspondence for $k= \bold F_q$
\endheading
Next we study the case $k= \bold F_q$.
\proclaim{Theorem 2}
Let $D=(R, \{l_i\}_i, \{\kappa_i\}_i)$ be a non-resonant 
hypergeometric data 
and $D_{red}$ be the hypergeometric data
defined in Proposition 5.1.  We renumber the index
such that $S=\{ 1, \dots ,s\}$.  Then there exists the following isomorphism
of constructible sheaves induced by an algebraic correspondence:
$$
\Cal G (D) \simeq \Cal G(D_{red}) \otimes 
H^s(FAS)(\kappa_1, \dots ,\kappa_s, \psi ).
$$
\endproclaim
\demo{Proof}
The proof is simpler compared to the case $k = \bold C$.
Let $y_i$ denote the element corresponding to $q(e_i)= \omega_i$
in $k[L]$ for $i=1, \dots ,s$. Take a base of 
$L_{red}\simeq \bold Z^{n-r}$. Then 
$k[L_{red}]$ is identified
with $k[x_{s+1}^{\pm}, \dots , x_n^{\pm}]$. 
We use a system of coordinate $\bar x=(x_{s+1}, \dots ,x_n)$.
If we put
$f_0 = \sum_{i=s+1}^r a_i \bar x^{\bar\omega_i}$, the variety $\Cal X_D$
is denoted by 
$$
\Cal X_D = \{ (y,\bar x,t ) \in \bold T(y) \times \bold T(x)
\times \bold A (t ) \mid 
a_1y_1+ \cdots + a_r y_r + f_0 =t \}
$$
We define the variety $\tilde\Cal X_D$ by
$$
\align
\tilde\Cal X_D  =  \{ &  (y,\bar x,t, s, a) \in
\bold T(y) \times \bold T(\bar x) \times \bold A(t) \times \bold A(s)
\times \bold T(a) \mid \\
&(a_1y_1 + \cdots + a_sy_s + f_0 (a,\bar x) = t ,
f_0 = s \} \\
\endalign
$$
It is easy to see that this is an etale covering of $\Cal X_D$ with
the Galois group $\bold F_p$.
We introduce a new coordinate $r =(r )$ and $Y=(Y_1, \dots ,Y_s)$.
By changing coordinate
$(y,\bar x, t , s ,a) \to (Y, \bar x, s , r ,a)$
given by $a_iy_i =Y_i$ and $r = t -s$, the variety
$\tilde\Cal X_D$ is isomorphic to $\Cal X_{D_{red}} \times FAS$, where
$$
\align
\Cal X_{D_{red}} = &\{ (x,\sigma , a) \in \bold T(x) \times \bold A(\sigma )
\times \bold T(a) \mid  f_0(a, x) = s \} \\
FAS = & \{ (Y,\rho ) \in \bold T(Y) \times \bold A (\tau ) \mid
Y_1 + \dots + Y_s = r \}. \\
\endalign
$$
The corresponding rank 1 local system is 
$pr_1^*(\Cal F^{(1)})\otimes pr_2^*(\Cal F^{(2)})$, where
$$
\Cal F^{(1)} = K (X_{D_{red}};\prod_{i=s+1}^r(a_ix^{\omega_i})^{\kappa_i})
\otimes \Cal L_{\psi}(s) \text{ and }
\Cal F^{(2)} = K (FAS ; \prod_{i=1}^sY_1^{\kappa_i})
\otimes\Cal L_{\psi}(r).
$$
Let $\varphi_{red}$ be the natural projection
$\Cal X_{D_{red}} \to \bold T(a_{s+1}, \dots ,a_r)$.
Since $\tilde\Cal X_{D_{red}}$ is etale over $\Cal X_{D_{red}}$,
we have
$$
\align
& R^n\varphi_*K(\Cal X_D;\prod_{i=1}^sY_i^{\kappa_i}
\prod_{i=s+1}^r(a_ix^{\omega_i})^{\kappa_i} )
\otimes \Cal L_{\psi}(t)  \\
&\simeq
R^{n-s}\varphi_{red,*}K(\Cal F^{(1)}) \otimes
H^{s}(FAS ; \kappa_1, \dots ,\kappa_s, \psi ). \\
\endalign
$$
Therefore we have the proposition.
\enddemo

\heading
\S 6 Algebraic correspondences associated to the multiplicative
equivalence
\endheading
\heading
\S 6.1 Existence of good base
\endheading

In this section we introduce an algebraic correspondence
associated to multiplicative equivalence.
Let $D= (R, \{l_i\}_i, \{ \kappa_i\}_i)$ be a hypergeometric data 
and suppose that $l_1(R) =d\bold Z$ with $d >0$.
\proclaim{Proposition 6.1}
There exists a base of $L\simeq \bold Z^{n+1}$ such that
$$
<\omega_i>_{i=2,\dots ,r} = 
<d\omega_1, \omega_i>_{i=2. \dots ,r}
= d\bold Z \oplus \bold Z^n
$$
and $\omega_1 = (k,0, \dots ,0)$, with $(k,d)=1$.
\endproclaim
\demo{Proof}
Let $Q = <\omega_i>_{i=2, \dots ,r}$. By the following commutative
diagram, we have $L /Q = \bold Z /d\bold Z$.
$$
\CD
@. 0 @. 0 @. 0 \\
@. @VVV @VVV @VVV \\
0 @>>> \bold Z^{r-1} \cap R @>>> \bold Z^{r-1} @>>> Q @>>> 0 \\
@. @VVV @VVV @VVV \\
0 @>>> R @>>> \bold Z^r @>>> L @>>> 0 \\
@. @VVV @VVV @VVV \\
0 @>>> \Im (l_1) @>>> \bold Z @. Q/L \\
@. @VVV @VVV @VVV \\
@. 0 @. 0 @. 0 \\
\endCD
$$
Here, columns and rows are exact.
If we put $Q' = Q+d\omega_1\bold Z$, then by the following commutative
diagram, we have $L/Q'\simeq \bold Z/d\bold Z$, and as consequence,
we have $Q = Q'$.
$$
\CD
0 @>>>  R @>>> d\bold Z \oplus \bold Z^{r-1} @>>> Q' @>>> 0 \\
@. @VVV @VVV @VVV \\
0 @>>> R @>>> \bold Z^r @>>> L @>>> 0 \\
\endCD.
$$
Let $\omega_1 = k l$ with a primitive element $l$ in $L$.
If we write $l= a_1\omega_1 + \dots a_r\omega_r$, then
$(1-ka_1)\omega_1 \in Q$, $(k, d)$=1.
The natural map $\bold Z \omega_1 \to L/Q$ is surjective,
the homomorphism $\bold Z l \to L/Q$ is also surjective. Therefore
$dl$ is a primitive element in $Q$, and we can take a base of $Q$;
$Q=\bold Z dl \oplus \oplus_{i=1}^n \bold Zu_i$.
\enddemo
We construct an algebraic correspondence using a coordinate
given in Proposition 6.1.
We a take base of $L$ : $L \simeq \bold Z^{n+1}$ which satisfies
the condition of Proposition 6.1.
Under this base, the first coordinate of $\omega_i$ is
divisible by $d$ for $i=2, \dots, r$. Therefore 
$f_0 = \sum_{i=2}^ra_i x^{\omega_i}$ can be written as
$f_0(x_1^d, x_2, \dots ,x_{n+1})$.  The polynomial $f$ can be written
as $a_1 x_1^k + f_0(x^d, x_2, \dots ,x_{n+1})$.
Using the base as in Proposition 6.2, the lattice 
$L^{(1,d)} = (\bold Z^d \oplus \bold Z^{r-1})/\bold l^{(1,d)}(R)$
and the morphism $\bold Z^d \oplus \bold Z^{r-1} \to L^{(1,d)}$
is expressed as follows.  We write the structure of $L^{(1,d)}$
multiplicatively. Let $u_1,\dots , u_d$ be a system coordinates.
$L^{(1,d)}$ is generated by $Q$ and $u_1,\dots ,u_d$ with
the relation $u_1 \cdots u_d = z^k$, where $z=x^l$ and
$l=(1, 0, \dots ,0)$. Note that $\alpha_1 = k\kappa_1$ under
this base. (For the definition of $\bold l^{(1,d)}$ and $L^{(1,d)}$
see Section 4.)
Before constructing algebraic correspondences, we note the following
lemma.
\proclaim{Lemma 6.2}
Let $d$ be a natural number prime to the characteristic of
$k$.
Let $\xi_1, \dots ,\xi_d$ be a set of variables and $s_1, \dots s_d$
be the elementary symmetric functions of $\xi_i$ ($i= 1, \dots, d$) 
of degree $1, \dots ,d$.  Then $\xi_1^d + \dots + \xi_d^d$
is expressed as a polynomial $F(s_1, \dots ,s_d)$ of $s_1, \dots ,s_d$
and we have
$$
F(s_1, \dots, s_d) = F(s_1, \dots ,s_{d-1} ,0) +(-1)^{d-1}d s_d.
$$  
Moreover, if $s_1, \dots , s_{d-1}$ are elementary symmetric
functions in $\gamma_1, \dots ,\gamma_{d-1}$, then
$$
F(s_1, \dots, s_{d-1}, 0) = \gamma_1^d +\dots + \gamma_{d-1}^d.
$$
\endproclaim
\demo{Proof}
See [T].
\enddemo
\heading
\S 6.2 Multiplicative correspondences for $k =\bold C$.
\endheading
The main theorem of this section is the following.
\proclaim{Theorem 3}
Let $D$ be a hypergeometric data with 
$l_i(R) = d\bold Z$ (if $k= \bold C$) or
$l_i(R) = dp^e \bold Z$ with $(d,p) =1$ (if $k=\bold F_q$) and $D^{(i,d)}$
be the hypergeometric data defined in Section 3. We put
$\tilde D^{(d)} = D \oplus \oplus_{i=1}^{d-1}(0,\frac{i}{d})$.
Then we have the following isomorphism of hypergeometric sheaves induced
by an algebraic correspondence.
\roster
\item
If $k=\bold C$,
$
\Cal G(\tilde D^{(d)}) \otimes K((\frac{(-1)^{d-1}}{d})^{\kappa_1})
\simeq T_{(i,d)}^* \Cal G(\tilde D^{(i,d)}),
$
where $K((\frac{(-1)^{d-1}}{d})^{\kappa_1})$ is the Hodge structure
on $\Spec \bold C$ corresponding to the Kummer character
$(\frac{(-1)^{d-1}}{d})^{\kappa_1}$.
\item
If $k= \bold F_q$,
$
\Cal G(\tilde D^{(d)},\psi ) \otimes K((\frac{(-1)^{d-1}}{d})^{\kappa_1})
\simeq T_{(i,d)}^* \Cal G(\tilde D^{(i,d)}, \psi ),
$
Here $T_{(i,d)}^*$ denotes the translation of $\bold T(R)$ given by
$T^*_{(i,d)}(a_i)=(-1)^{d-1}a_1/d$ and
$K((\frac{(-1)^{d-1}}{d})^{\kappa_1})$ is the $l$-adic sheaf on
$\Spec \bold F_q$ corresponding to the Kummer character
$(\frac{(-1)^{d-1}}{d})^{\kappa_1}$.
\endroster
\endproclaim
\demo{Proof}
In this subsection, we prove (1).
We choose a base of $L$ satisfying the condition of Proposition 6.1.
Let $\kappa = (\kappa_1, \dots ,\kappa_r)$ and 
$\alpha =(\alpha_1, \dots, \alpha_{n+1}) \in L \otimes I(\bold C, K)$
be the image of $\kappa$
under the homomorphism $q$.
We prove the above theorem for the case 
$\sum_{i=1}^r \kappa \neq 0$.  In the case 
$\sum_{i=1}^r \kappa = 0$, we can prove the theorem by similar argument.
We introduce systems of variables 
$$
\align
& x_1,
x = (x_2, \dots ,x_{n+1}),
\xi = (\xi_1, \dots , \xi_d),
\gamma = (\gamma_1, \dots, \gamma_{d-1}), \\
& u=(u_1, \dots ,u_d), g =(g_1,\dots, g_{d-1}),
s^0 = (s_1, \dots ,s_{d-1}),s=(s_1, \dots ,s_d), z \\
& a=(a_1,\dots ,a_r), a_{\geq 2} = (a_2, \dots ,a_r),
b=(b_1, \dots ,b_d), c =(c_1, \dots ,c_{d-1}) \\
\endalign
$$
We use the notation $a'_1 =(-1)^{d-1}a_1/d$ for short.
We define five varieties $\Cal X_{D^{(1,d)}}$, $\Cal X_1$, 
$\Cal Y$, $\Cal X_2$ and $\Cal X_{\tilde D^{(d)}}$.
$$
\align
\Cal X_{D^{(1,d)}} = & \{ (u , x,z,b, a_{\geq 2} ) \in \bold T(u) 
\times \bold T(x) \times \bold T(z) \times
\bold T(b) \times \bold T(a_{\geq 2}) \mid \\
& b_1u_1 + \cdots + b_du_d + 
f_0(z, x_2, \dots ,x_{n+1}) = 1,u_1\cdots u_d = z^k \} \\
\Cal X_1 = &\{ (\xi , x,z, a ) \in \bold T(\xi) \times \bold T(x) \times
\bold T(z) \times \bold T(a) \mid \\
& a_1 (\xi_1^d + \cdots + \xi_d^d) + 
f_0(z, x_2, \dots ,x_{n+1}) = 1,
(\xi_1\cdots \xi_d)^d =z^k\} \\
\Cal Y = &\{ (s^0, s_d , x,z, a ) \in \bold A(s^0) \times 
\bold T(s_d)\times \bold T(x) \times \bold T(z)\times
\bold T(a) \mid \\
& a'_1 (F(s_1, \dots ,s_{d-1},0) + (-1)^{d-1}ds_d) + 
f_0(z, x_2, \dots ,x_{n+1}) = 1,
s_d^d = z^k \} \\
\Cal X_2 = &\{ (\gamma ,x_1, x, a ) \in \bold A(\gamma ) \times 
\bold T (x_1) \times \bold T(x) \times
\bold T(a) \mid \\
& a'_1 (\gamma_1^d + \cdots + \gamma_{d-1}^d) + 
a_1 x_1^k + f_0(x_1^d, x_2, \dots ,x_{n+1}) = 1 \} \\
\Cal X_{\tilde D^{(d)}} = &\{ (g ,x_1, x, c, a ) \in \bold T(g ) \times 
\bold T (x_1) \times \bold T(x) \times
\bold T(c) \times\bold T(a) \mid \\
& c_1g_1 + \cdots + c_{d-1}g_{d-1} + 
a_1 x_1^k + f_0(x_1^d, x_2, \dots ,x_{n+1}) = 1 \} \\
\endalign
$$
Let $(\xi , x,z, a )$ (resp. $(\gamma,x_1, x, a)$) be a point in $\Cal X_1$
(resp. $\Cal X_2$).
By substituting
$s_i$ ($i=1, \dots ,d$ (resp. $i=1, \dots ,d-1$))
by the $i$-th symmetric polynomial
of $\xi_i$ ($i=1, \dots ,d$) (resp. $\gamma_i$ (
$i=1. \dots , d-1$)),
we get a finite homomorphism $\pi_1$ (resp. $\pi_2$)
from $T_{(1,d)}^*(\Cal X_1)$ to $\Cal Y$ (resp. $\Cal X_2$ to $\Cal Y$)
using Lemma 6.2. Then $\pi_1$ and $\pi_2$ are
quotient by finite groups $\frak S_d$ and $\frak S_{d-1}$ respectively.
The natural morphisms $\Cal X_1 \to \bold T(a)$, 
$\Cal Y \to \bold T(a)$ and $\Cal X_2 \to \bold T(a)$ are
denoted by $\phi_1$, $\phi_y$ and $\phi_2$ respectively.
We consider the following sheaves on each varieties.
$$
\align
\Cal F^{(1,d)} & =K (\Cal X_{(1,d)}; 
\prod_{i=1}^d(b_iu_i)^{(i+\kappa_1)/d}
\prod_{i=2}^{n+1}x_i^{\alpha_i}\prod_{i=2}^r a_i^{\kappa_i} ) \\
\Cal F_1 & = K (\Cal X_1; 
\prod_{i=1}^d (a_1\xi_i^d)^{(i+\kappa_1)/d}
\prod_{i=2}^{n+1}x_i^{\alpha_i}\prod_{i=2}^r a^{\kappa_i} ) \\
\Cal F_y & = K (\Cal Y;(a')^{\frac{d-1}{2}} 
(a_1s_d)^{\kappa_1}
\prod_{i=2}^{n+1}x_i^{\alpha_i}\prod_{i=2}^r a^{\kappa_i} ) \\
\Cal F_2 & = K (\Cal X_2; \prod_{i=1}^{d-1}(a_1'\gamma_i^d )^{i/d}
(a_1x_1^k)^{\kappa_1} 
\prod_{i=2}^{n+1}x_i^{\alpha_i}\prod_{i=2}^r a_i^{\kappa_i} ) \\
\Cal F^{(d)} & = \Cal F (\Cal X_{\tilde D^{(d)}}; 
\prod_{i=1}^{d-1}(c_ig_i)^{i/d}(a_1 x_1^k)^{\kappa_1} 
\prod_{i=2}^{n+1}x_i^{\alpha_i}\prod_{i=2}^r a_i^{\kappa_i} ) \\
\endalign
$$
By the argument in \S 2, the variety 
$\Cal X_{D^{(1,d)}} \to \bold T(b) \times \bold T(a_{\geq 2})$ and
$\Cal X_{\tilde D^{(d)}} \to \bold T(c) \times \bold T(a)$ descend
to variety $\Cal X_R \to \bold T(R)$ and $\Cal X_R' \to \bold T(R)$
respectively.  Under this descent, the sheaves $\Cal F^{(1,d)}$
and $\Cal F^{(d)}$ descend to sheaves $\bar\Cal F^{(1,d)}$
and $\bar\Cal F^{(d)}$ on $\Cal X_R$ and $\Cal X'_R$ respectively.  
By the definition of hypergeometric sheaves, we have
$$
\align
R^n\varphi_{R*}\bar\Cal F^{(1,d)} & = \Cal G(D^{(1,d)}) \\
R^n\varphi'_{R*}\bar\Cal F^{(d)} & = \Cal G(\tilde D^{(d)}). \\
\endalign
$$
We consider the following diagram:
$$
\CD
\Cal X_1 @>>> \Cal X_{D^{(1,d)}} @>>> \Cal X_R \\
@V{\varphi_1}VV  @VVV @V{\varphi_R}VV \\
\bold T(a) @>{\sum}>> \bold T(b)\times \bold T(a_{\geq 2}) @>>> \bold T(R) \\
(a_1, \dots ,a_r) &\mapsto &(b_1, \dots, b_r, a_2, \dots ,a_r)  \\
& &  =(a_1, \dots, a_1, a_2, \dots ,a_r) \\
\endCD
$$
$$
\CD
\Cal X_2 @>>> \Cal X_{\tilde D^{(d)}} @>>> \Cal X'_R \\
@V{\varphi_2}VV  @VVV @V{\varphi'_R}VV \\
\bold T(a) @>{\sum '}>> \bold T(c)\times \bold T(a) @>>> \bold T(R), \\
(a_1, \dots ,a_r) &\mapsto &(c_1, \dots, c_{r-1}, a_1, \dots ,a_r) \\
& &  =(a_1, \dots, a_1, a_1, \dots ,a_r) \\
\endCD
$$
where $\Cal X_1 \to \Cal X_{D^{(1,d)}}$ and 
$\Cal X_2 \to \Cal X_{\tilde D^{(d)}}$ are given by
$$
(u,x,z,b,a_{\geq 2})=
((\xi_1^d,\dots ,\xi_d^d),x,z,(a_1, \dots ,a_1),a_{\geq 2})
$$
and
$$
(g,x_1,x,c,a) = ((\gamma_1^d,\dots ,\gamma_{d-1}^d),x_1,x,(a_1, \dots ,a_1),a).
$$
Let $\bold l^* :\bold T(a) \to \bold T(R)$ be the morphism induced by $\bold l$.
Since
$(\sum^* \Cal X_{D^{(1,d)}})\times_{\bold T(u)}\bold T(\xi ) 
= \Cal X_1$ and 
$((\sum ')^* \Cal X_{\tilde D^{(d)}})\times_{\bold T(g)}\bold T(\gamma ) =
\Cal X_2$,
we have
$$
\align
\bold l^* \Cal G(D^{(1,d)}) & \simeq R^n \varphi_{1*}\Cal F_1 \\
\bold l^* \Cal G(\tilde D^{(d)}) & \simeq R^n \varphi_{2*}\Cal F_2. \\
\endalign
$$
Therefore to prove the theorem,
it is enough to show the following theorem.
\enddemo
\proclaim{Proposition 6.3}
On the torus $\bold T(a)$, the morphism $\pi_1$ and $\pi_2$ 
induce isomorphisms of sheaves:
$$
\align
T^*_{(1,d)}(R^n\varphi_{1*}\Cal F_1)  & \simeq
R^n\varphi_{y}\Cal F_{y} \otimes
K((\frac{(-1)^{d-1}}{d})^{\kappa_1})
\text{ and } \tag{1}\\
R^n\varphi_{1*}\Cal F_2  & \simeq
R^n\varphi_{y}\Cal F_{y}  \tag{2}\\
\endalign
$$
Moreover they are compatible with the descent data.
\endproclaim
\demo{Proof}
We prove the isomorphism (1).  The isomorphism (2) can be proved
similarly. 
Let 
$$
B=\{(T,D,x,a) \in \bold A(T)\times \bold T(D) \times \bold T(x) \times
\bold T(a) \mid  -a_1'T + f_0(D, x_2, \dots ,x_{n+1}) =1\}
$$
be a variety over $\bold T(a)$.  The natural map $B \to \bold T(a)$
is denoted by $\varphi_B$.
Then the morphism 
$$
(\xi, x,z,a) \mapsto 
(\xi_1^d+ \dots + \xi_d^d,z,x,a)
\text{ and }
(s^0,s_d, x,z,a) \mapsto (F(s_0, \dots ,s_d),z,x,a)
$$
induces a morphism
from $h'_1 :T^*_{(1,d)}\Cal X_1 \to B$ and 
$h_y : \Cal Y \to B$.  
Consider the following commutative diagram;
$$
\CD
T^*_{(1,d)}\Cal X_1 @>{\pi_1}>> \Cal Y \\
@V{h'_1}VV @VV{h_y}V \\
B @= B @>>{\varphi_B}> \bold T(a) \\
\endCD
$$
Let $\Cal F_B$ be
a sheaf on $B$ defined by 
$$
\Cal F_B = K ( B; (a'_1)^{\frac{d-1}{2}}(a_1D^{\frac{k}{d}})^{\kappa_1}
\prod_{i=2}^nx_i^{\alpha_i}\prod_{i=2}^ra_i^{\kappa_i}).
$$
We have
$$
\align
\bold R\varphi_{B!}(\bold R h'_{1!}K(1,\dots ,d)\otimes \Cal F_B)
\otimes K((\frac{(-1)^{d-1}}{d})^{\alpha_1}) 
& \simeq T^*_{(1,d)}\bold R \varphi_{1!}\Cal F_1 
\tag{1}\\
\bold R \varphi_{B!}(\bold R h_{y!}K\otimes \Cal F_B)
& \simeq \bold R \varphi_{y!}\Cal F_y. 
\tag{2}\\
\endalign
$$
By the Poincare duality, it is enough to prove that the 
morphism $\pi_1$ induces an isomorphism between (1) and (2).
By the comparison theorem, it is enough to prove the isomorphism
for the etale cohomology
with the coefficient $\bar\bold Q_l$ and by the specialization
argument, it is enough to show the isomorphism by reduction 
mod $p$. Suppose that all the varieties, sheaves and
automorphisms
are defined over a finite filed $\bold F_q$ and $d \mid (q-1)$.
Since $T^*_{(1,d)}\Cal X_1 \to \Cal Y$ is a finite $\frak S_{d}$
covering, we have
$
\bold Rh_{y!}K \simeq (\bold R h'_{1!}K)^{\frak S_d} 
$
and the action of $(\mu_d)^d$ on $\Cal X_1$ induces a decomposition
$
\bold R h'_{1!} K \simeq \oplus_{(t_1,\dots ,t_d)}
\bold R h'_{1!} K (t_1, \dots ,t_d)
$
according to the characters $(t_1, \dots ,t_d)$ of $(\mu_d)^d$.
For an element $\sigma \in \frak S_d$, we have
$\sigma^* (\bold R h'_{1!} K (t_1, \dots ,t_d)) =
\bold R h'_{1!} K (t_{\sigma (1)}, \dots ,t_{\sigma (d)})$.
\enddemo
\proclaim{Lemma 6.4}
Suppose that $t_p = t_q$ for some $1 \leq p \neq q \leq d$.
Let $\sigma_{p,q}$ be the transposition of $p$ and $q$.
Then for a closed point $p=(T,D)$ of $B$, we have
$$
\tr (Frob^m_{\kappa (p)}(\sigma_{p,q}+1)
\mid \bold R h'_{1!}K(t_1, \dots ,t_d))_{\bar\kappa (p)}=0
$$
for all $m \geq 1$.
\endproclaim
The proof of this lemma is given later.  
Let 
$
\Cal F'_1(t_1, \dots ,t_d) = 
\bold R h'_{1!}K(t_1,\dots ,t_d). 
$
If $t_p = t_q$, by Lemma 6.4 and
Lefschetz trace formula for etale cohomology, we have
$$
\align
&
\tr (Frob^m_{\kappa (a)}(\sigma_{p,q}+1)
\mid \bold R \varphi_{B!}\Cal F'_1(t_1, \dots ,t_d))_{\bar\kappa (a)} \\
= &
\sum_{p \in \varphi_B^{-1}(a)(\kappa^m (a))}
(\tr (Frob_p(\sigma_{p,q}+1)
\mid \bold R h'_{1!}K(t_1, \dots ,t_d))_{\bar\kappa (p)})
\cdot Frob_p\mid (\Cal F_B)_{\bar\kappa (p)}  \\
= & 0. \\
\endalign
$$
Here $\kappa^m(a)$ denotes the degree $m$ extension of $\kappa (a)$.
Therefore by Proposition 2.4, we have
$$
\tr (Frob^m_{\kappa (a)}(\sigma_{p,q}+1)
\mid R^n \varphi_{B!}\Cal F'_1(t_1, \dots ,t_d))_{\bar\kappa (a)} =0
$$
for all positive $m$ and $\sigma_{p,q}$ acts as $-1$ on
$R^n \varphi_{B!}\Cal F'_1(t_1, \dots ,t_d))_{\bar\kappa (a)}$.
As a consequence, we have
$$
(R^n \varphi_{B!}\bold R h'_{1!}K)_{\bar\kappa (a)}^{\frak S_d} \simeq
R^n \varphi_{B!}\Cal F'_1(1, \dots ,d))_{\bar\kappa (a)}
$$
and we obtain the proposition.
\demo{Proof of Lemma 6.4}
We may assume that $\kappa (p) = \bold F_q$.
Let $X_{T,D_0}$ be the variety defined by
$$
X_{T,D_0} = \{ \xi \in \bold T(\xi )  \mid
\xi_1^d + \cdots + \xi_d^d = T, (\xi_1\cdots \xi_d)^d =D_0 \}.
$$
It is easy to see that the fiber of $h'_1$ at $(T,D)$ is
isomorphic to $X_{T,D^k}$.
Though this variety is not geometrically connected, we consider the 
cohomology with compact support.  For a character
$t=(t_1, \dots, t_d)$ of $(\mu_d)^d$, 
$H^i_c(\bar X_{T,D}, K)(t_1, \dots ,t_d)$,
denotes the $t$-part of the corresponding cohomology.
Define $\Phi (T)$ by 
$$
\Phi (T) =\sum_{i=0}^{2d-4}(-1)^i\tr (Frob_{\bold F_q}(\sigma_{p,q} +1)
\mid
H^i_c(\bar X_{T,D}, K)(t_1, \dots ,t_d)).
$$
First we assume that $T \neq 0$.
Let $\Cal Z_{D_O}$ be the variety defined by
$$
\Cal Z_{D_0} = \{ (\xi, T) \in \bold T(\xi) \times \bold T(T) \mid
\xi_1^d +\cdots + \xi_d^d = T, (\xi_1 \cdots \xi_d)^d = D_0 \}
$$
and $u : \Cal Z_{D_0} \to \bold T(T)$ be the natural homomorphism.
For a multiplicative character $\gamma : \bold F_q^{\times} \to K$,
$K(\gamma )$ denotes the corresponding rank 1 local system
on $\bold T(T)$.  Then by the Lefschetz trace formula,
we have
$$
\tr (Frob_{\bold F_q}(\sigma_{p,q} +1) \mid
H^*_c(\bar\bold T(T), \bold R u_! K(t_1, \dots ,t_d) \otimes K(\gamma ))
=\sum_{T \in \bold F_{q}^{\times}}\Phi(T)\gamma (T).
$$
On the other hand, 
$$
\align
& H^i_c(\bar\bold T(T), \bold R u_! K(t_1, \dots ,t_d) \otimes K(\gamma )) \\
= &
H^i_c(\Cal Z_{D_0}, K)(t_1, \dots, t_d,\gamma ) \\
=& \cases
\text{ $1$-dimensional vector space over $K$ } & \qquad 
(\text{ if } i = d-1 ) \\
0 & \qquad
(\text{ if } i \neq d-1 ) \\
\endcases
\endalign
$$
since $\Cal Z_{D_0}$ is a quotient of a $d-1$-dimensional Fermat hypersurface.
It is easy to see that $\sigma_{p,q}$ acts as $(-1)$-multiplication
on $H^{d-1}(\Cal Z_{D_0}, K)(t_1, \dots, t_d,\gamma )$.
Therefore we have
$\sum_{T \in \bold F_{q}^{\times}}\Phi(T)\gamma (T)=0$ for all $\gamma$.
As a consequence, we have $\Phi (T) =0$.  For the case, $T=0$,
the variety $X_{0,D}$ is a quotient of a Fermat hypersurface.
In this case we also have $\Phi (0) =0$.
\enddemo

\heading 
\S 6.3 The case $k=\bold F_q$
\endheading
In this subsection, we prove (2) of Theorem 3.
We use the same systems of coordinates. 
Let $\# \bold Z/\Im (R)=dp^e$, where $(d,p) =1$.
As in \S 6.2, we define 
five varieties $\Cal X_{D^{(1,d)}}$, $\Cal X_1$,
$\Cal Y$, $\Cal X_2$ and $\Cal X_{\tilde D^{(d)}}$ by
$$
\align
\Cal X_{D^{(1,d)}} = & \{ (u , x,b,z,t, a_{\geq 2} ) \in \bold T(u) 
\times \bold T(x) \times
\bold T(b) \times \bold T(z) \times \bold A(t) \times
\bold T(a_{\geq 2}) \mid \\
& b_1u_1 + \cdots + b_du_d + 
f_0(t^{p^e}, x_2, \dots ,x_{n+1}) = t,
u_1 \cdots u_d = z^k \} \\
\Cal X_1 = &\{ (\xi , x,z,t, a ) \in \bold T(\xi) \times \bold T(x) \times
\bold T(z) \times \bold A(t) \times \bold T(a) \mid \\
& a_1 (\xi_1^d + \cdots + \xi_d^d) + 
f_0(t^{p^e}, x_2, \dots ,x_{n+1}) = t,
u_1 \cdots u_d = z^k \} \\
\Cal Y = &\{ (s^0, s_d , x,z,t, a ) \in \bold A(s^0) \times 
\bold T(s_d)\times \bold T(x) \times \bold T(z)
\times \bold A(t) \times \bold T(a) \mid \\
& a'_1 (F(s_1, \dots ,s_{d-1},0) + (-1)^{d-1}ds_d) + 
f_0(z^{p^e}, x_2, \dots ,x_{n+1}) = t,
s_d^d = t^k \} \\
\Cal X_2 = &\{ (\gamma ,x_1, x, t,a ) \in \bold A(\gamma ) \times 
\bold T (x_1) \times \bold T(x) \times
\bold A(t) \times\bold T(a) \mid \\
& a'_1 (\gamma_1^d + \cdots + \gamma_{d-1}^d) + 
a_1 x_1^k + f_0(x_1^{dp^e}, x_2, \dots ,x_{n+1}) = t \} \\
\Cal X_{\tilde D^{(d)}} = &\{ (g ,x_1, x,t, c, a ) \in \bold T(g ) \times 
\bold T (x_1) \times \bold T(x) \times \bold A(t)\times
\bold T(c) \times\bold T(a) \mid \\
& c_1g_1 + \cdots + c_{d-1}g_{d-1} + 
a_1 x_1^k + f_0(x_1^{dp^e}, x_2, \dots ,x_{n+1}) = t \} \\
\endalign
$$
We define a morphism $\pi_1: T_{(1,d)}^* \Cal X_1 \to \Cal Y$
and $\pi_2 : \Cal X_2 \to \Cal Y$ as in \S 6.2.  Then $\pi_1$
and $\pi_2$ are quotient by the symmetric group $\frak S_d$
and $\frak S_{d-1}$.  We define sheaves on each varieties:
$$
\align
\Cal F^{(1,d)} & =\kappa (\Cal X_{(1,d)}; 
\prod_{i=1}^d(b_iu_i)^{(i+\kappa_1)/d}
\prod_{i=2}^{n+1}x_i^{\alpha_i}
\prod_{i=2}^r a_i^{\kappa_i} )\otimes \Cal L_{\psi}(t) 
\\
\Cal F_1 & = \kappa (\Cal X_1; 
\prod_{i=1}^d (a_1\xi_i^d)^{(i+\kappa_1)/d}
\prod_{i=2}^{n+1}x_i^{\alpha_i}
\prod_{i=2}^r a^{\kappa_i} )\otimes \Cal L_{\psi}(t) 
 \\
\Cal F_y & = \kappa (\Cal Y;(a')^{\frac{d-1}{2}} 
(a_1s_d)^{\kappa_1}
\prod_{i=2}^{n+1}x_i^{\alpha_i}
\prod_{i=2}^r a^{\kappa_i} )\otimes \Cal L_{\psi}(t) 
 \\
\Cal F_2 & = \kappa (\Cal X_2; \prod_{i=1}^{d-1}(a_1'\gamma_i^d )^{i/d}
(a_1x_1^k)^{\kappa_1} 
\prod_{i=2}^{n+1}x_i^{\alpha_i}
\prod_{i=2}^r a_i^{\kappa_i} )\otimes \Cal L_{\psi}(t) 
 \\
\Cal F^{(d)} & = \kappa (\Cal X_{\tilde D^{(d)}}; 
\prod_{i=1}^{d-1}(c_ig_i)^{i/d}(a_1 x_1^k)^{\kappa_1} 
\prod_{i=2}^{n+1}x_i^{\alpha_i}
\prod_{i=2}^r a_i^{\kappa_i} )\otimes \Cal L_{\psi}(t) 
 \\
\endalign
$$
Let $\bold l^*$ be the morphism $\bold T (a) \to \bold T(R)$
induced by $\bold l$.  Then we have
$$
\align
\bold l^* \Cal G(D^{(1,d)},\psi ) & \simeq R^n \varphi_{1*}\Cal F_1 \\
\bold l^* \Cal G(\tilde D^{(d)}, \psi ) & \simeq R^n \varphi_{2*}\Cal F_2. \\
\endalign
$$
As in \S 6.2, it is enough to prove the following proposition.
\proclaim{Proposition 6.5}
On the torus $\bold T(a)$, the morphism $\pi_1$ and $\pi_2$ 
induce isomorphisms of sheaves:
$$
\align
T^*_{(1,d)}(R^n\varphi_{1*}\Cal F_1)  & \simeq
R^n\varphi_{y}\Cal F_{y}\otimes K((\frac{(-1)^{d-1}}{D})^{\alpha_1})
 \text{ and } \tag{1}\\
R^n\varphi_{1*}\Cal F_2  & \simeq
R^n\varphi_{y}\Cal F_{y}  \tag{2}\\
\endalign
$$
Moreover they are compatible with the descent data.
\endproclaim
\demo{Proof}
The proof is similar to that of Proposition 6.3 and omit the proof.
\enddemo

\heading
\S 7 Multiplication by $p$ for $k= \bold F_q$
\endheading
Let $D=(R, \{ l_i\}_i, \{\kappa_i\}_i)$ be a hypergeometric data
and we put $l_1(R)=d\bold Z$.  Suppose that $d$ is divisible by $p$;
$d = pd'$. Since $p$ is invertible in $I(\bar\bold F_q, \bar\bold Q_l)$,
there exists unique $\kappa'_i$ such that $p\kappa'_i = \kappa_i$.
Let $l_1' =\frac{1}{p}l_1$ then the non-resonance condition for 
$D=(R, \{l_i\}_i, \{\kappa_i\}_i)$ is equivalent to that for
$D^{(1,p)} = (R, \{l'_1,l_i\}_{i\geq 2},\{\kappa'_1,\kappa_i\}_{i\geq 2})$
by the following commutative diagram.
$$
\CD
R @>{(l_1',l_2, \dots ,l_r)}>> \bold Z \oplus \bold Z^{r-1} @>>> L' \\
@VVV @V{(p,1, \dots ,1)}VV @VVV   \\
R @>>{(l_1, \dots ,l_r)}> \bold Z \oplus \bold Z^{r-1} @>>> L  \\
\endCD
$$
By choosing a sufficiently good base as in the last section,
we may assume that $\omega_1 = (k,0 \dots , 0)$ and
the first component of $\omega_i$ ($i=2, \dots ,n$) is divisible
by $d=pd'$. Therefore we can write $\sum_{i=2}^r a_ix^{\omega_i}
=f_0(x_1^p,x_2, \dots ,x_{n+1})$ and 
the defining equation of $\Cal X$ and $\Cal X'$ for the hypergeometric
data $D$ and $D'$ is
$$
\align
\Cal X =&\{(x, \tau ,a) \in \bold T(x) \times \bold A (\tau )
\times \bold T(a) \mid
a_1 x_1^k + f_0(x_1^p, x_2, \dots ,x_{n+1})=\tau^p -\tau \} \\
\Cal X' =&\{(y,x_2, \dots x_n, \tau ,b_1, a_2, \dots ,a_r)
 \in \bold T(x) \times \bold A (\tau ) \times \bold T(b_1)
\times \bold T(a_2, \dots ,a_r) \\
& \mid
b_1 y_1^k + f_0(y_1, x_2, \dots ,x_{n+1}) =\sigma^p -\sigma \}. \\
\endalign
$$
The corresponding rank 1 local system is given by 
$$
K( \prod_{i=1}^r(a_ix^{\omega_i})^{\kappa_i})\text{ and } 
K((b_1y_1^k)^{\kappa'_1}
\prod_{i=2}^r(a_ix^{\omega_i})^{\kappa_i}). 
$$
Define a homomorphism $\Cal X \to \Cal X'$ by sending
$$
(x, \tau, a) \to (y_1, x_2, \dots ,x_{n+1}, \sigma , a) =
(x_1^p, x_2, \dots ,x_{n+1}, \tau + a_1x_1^k, a_1^p, a_2, \dots ,a_r).
$$
This is equivariant under the action of $\bold F_q$
Then it easy to see that the corresponding character coincides.
Taking account into the above commutative diagram and the definition
of descent data,
we have an isomorphism of hypergeometric sheaves on $\bold T(R)$:
$
\Cal G(D) \simeq \Cal G(D^{(1,p)})
$
by taking the Artin-Shreier character part of the higher direct images.
As a consequence, we have the following proposition.
\proclaim{Proposition 7.1}
Let $D=(R, \{ l_i\}_i, \{\kappa_i\}_i)$ be a 
non-resonant hypergeometric data
and such that $l_1(R)$ is divisible by $p$.  If we put
$D^{(1,p)} = (R, \{l'_1,l_i\}_{i\geq 2},\{\kappa'_1,\kappa_i\}_{i\geq 2})$,
where
$p\kappa'_i = \kappa_i$ and $l_1' =\frac{1}{p}l_1$.
Then $D^{(1,p)}$ is also non-resonant and on $\bold T(R)$,
we have an isomorphism of constructible sheaves
$
\Cal G(D) \simeq \Cal G(D^{(1,p)})
$
induced by an algebraic correspondences.
\endproclaim
\demo{Definition} The equivalence relation generated by $D\sim D^{(1,p)}$
is called the Frobenius equivalence.
\enddemo

\heading
\S 8 Cohomological Mellin Transform
\endheading
In this section, we give a definition of cohomological Mellin transform.
This notion is also introduce in [L-S] and [G-L] in different
presentation. First we consider the case $k= \bold C$.
\proclaim{Proposition 8.1}
Let $\Cal G (D)$ be the hypergeometric sheaf of a hypergeometric data
$D$ and $\chi \in R \otimes I(\bold C, K)$ be a character of 
$\pi_1(\bold T(R))$ of finite order.  Then the cohomology group 
$H^{r-n-1}(\bold T(R), \Cal G(D) \otimes K (\chi ))$ has a mixed
Hodge structure.  Moreover if $\kappa_i + l_i(\chi )$ are not
zero for all $i = 1,\dots ,r$, it is pure of weight $r$ (resp. $r-1$)
and one dimensional. Moreover the Hodge type is given by
$$
\align
& (\sum_{i=1}^r<\sigma(\kappa_i + l_i(\chi ))> + 
<-\sum_{i=1}^r \sigma\kappa_i>,
\sum_{i=1}^r<-\sigma(\kappa_i + l_i(\chi ))> + 
<\sum_{i=1}^r \sigma\kappa_i>) \\
& (resp. (\sum_{i=1}^r<\sigma(\kappa_i + l_i(\chi ))>,
\sum_{i=1}^r<-\sigma (\kappa_i + l_i(\chi ))>))\\
\endalign
$$
provided that $\sum_{i=1}^r \kappa_i \neq 0$ (resp. 
$\sum_{i=1}^r \kappa_i = 0$).
\endproclaim
\demo{Definition} 
For a constructible sheaf $\Cal F$
which is a mixed Hodge complex
on a torus $\bold T$ and an
element $\chi$ in $\Hom ( \bold T, \bold G_m)\otimes I(\bold C, K)$,
the mixed Hodge structure
\linebreak
$H^{r-n-1}(\bold T(R), \Cal F \otimes K (\chi ))$ is
denoted by $C(\Cal F, \chi )$.  This ``function `` from 
\linebreak
$\Hom ( \bold T, \bold G_m) \otimes I(\bold C, K)$ 
to the category of mixed Hodge structure
is called the cohomological Mellin transform of $\Cal F$.
\enddemo
\demo{Proof}
We choose a base of $L \simeq \oplus_{i=0}^{n}\bold Z v_i$.
We consider the variety over $\bold T(a)$.  By changing coordinate,
$\Cal X$ is actually defined over $\bold T(R)$.  To study the descent 
variety $\Cal X_R$, we choose a section $s : L \to \bold Z^r$ of 
$q : \bold Z^r \to L$.  Let $u_i =a^{s(v_i)}$ ($i=0,\dots,n$) and put
$y_i = u_ix_i$.  Then the equation for $\Cal X_R$ is 
$\sum a^{e_i-(s\circ q)(e_i)}y^{\omega_i}=1$.  Therefore the coefficient
$a^{e_i-(s\circ q)(e_i)} \in \bold C[R]$.  Using new coordinate 
$y= (y_0, \dots ,y_n)$, 
$$
\Cal X_R = \{ (r, y) \in \bold T(R) \times \bold T(y) \mid
\sum a^{e_i-(s\circ q)(e_i)}y^{\omega_i}=1 \}.
$$
Again we change coordinate given by 
$b_i = a^{e_i-(s\circ q)(e_i)}y^{\omega_i}$.  This change of coordinate
gives an isomorphism $\bold T(R) \times \bold T(y) \simeq \bold T(b)$,
where $b=(b_0, \dots ,b_n)$.  In fact,  we have an equality
$$
\cases
b^{e_i-(s\circ q)(e_i)} = a^{e_i-(s\circ q)(e_i)}\\
b^{-s(v_i)} = y_i. \\
\endcases
$$
Therefore $\Cal X_R \simeq \Cal X_b = \{ b \in \bold T(b) \mid
\sum_{i=1}^rb_i = 1\}$.  Under this isomorphism, the character 
$\kappa (\Cal X_R, x^\alpha a^\kappa )$ corresponds to
$\kappa (\Cal X_b, b^\kappa )$. We compute the cohomology
$H^{r-1}(\Cal X_b,\Cal F_{\chi} )$
using Leray spectral sequence for 
$\bar f : \Cal X_b \simeq \Cal X_R \to \bold T(R)$,
where $\Cal F_{\chi} = K(\kappa (\Cal X_b; b^{\kappa + \bold l (\chi )}))$.
By Proposition 2.5, $R^i\bar f_* \Cal F_{\chi}=0$ if $i\neq n$. Therefore
we have
$$
\align
H^{r-1} (\Cal X_b, \Cal F_{\chi}) \simeq &
H^{n-r-1}(\bold T(R), R^n\bar f_* \Cal F_{\chi}) \\
&
H^{n-r-1}(\bold T(R), \Cal G(D) \otimes K(\chi )). \\
\endalign
$$
To show the last part of the theorem, we can compute the Hodge number
of the cohomology of $H^{r-1}(\Cal X_b, \Cal F_{\chi})$ via that of
Fermat hypersurfaces (See [S].).
\enddemo

To state the similar results for $k=\bold F_q$, we recall the definition
of Gaussian sum.  Let $\bar p$ be an extension of valuation of $p$
to $\bar\bold Q \subset \bar\bold Q_l$.   We normalize additive valuation by 
$\ord_{\bar p} (p) =1$.  Via the reduction mod $\bar p$, 
$\mu_m(\bar\bold Q_l)$ is naturally identifies with $\mu_m(\bar\bold F_q)$
if $(m,p) =1$.  
The inverse of the reduction map is denoted by
$\omega : \bar\bold F_q \to \cup_{(m,p)=1}\mu_m(\bar\bold Q_l)$.
Therefore $I_m (\bar\bold F_q, \bar\bold Q_l)$ is identified
with $\frac{1}{m}\bold Z / \bold Z$, and $q=p^e$ such that $m \mid q-1$,
we define a Gaussian sum by 
$$
g(\kappa , \psi )=
g(\bold F_q , \kappa , \psi ) = \sum_{x \in {\bold F_q}^{\times}}
\chi_{\kappa} (x^{\frac{q-1}{m}})\psi (\tr_{\bold F_q/\bold F_p}(x)),
$$
where $\chi_{\kappa}(x) = \omega (x^{\frac{(q-1)\kappa}{m}})$.
It is easy to see that $g(0, \psi ) =-1$. 
The order of $g(\bold F_q, \kappa , \psi )$ at $\bar p$
is given by
$\ord_{\bar p} ( g(\bold F_q , \kappa , \psi )) =
\sum_{i=0}^{e-1}<p^i\kappa >$. Therefore 
$\frac{1}{e} \sum_{i=0}^{e-1}<p^i\kappa >$ is independent of the choice of
$q=p^e$ such that $m \mid q-1$.  Moreover if $\kappa = 0$, we have
$\ord_{\bar p} ( g(\bold F_q ,0 , \psi ) = 0$.
For $\sigma \in Gal (\bar\bold Q /\bold Q)$
$$
\ord_{\bar p} ( g(\bold F_q , \sigma\kappa , \psi )) =
\sum_{i=0}^{e-1}<tp^i\kappa >,
$$
where $t \in (\bold Z / m \bold Z)^{\times}$ is defined by the equality
$\sigma (\zeta ) = \zeta^t$ for all $\zeta \in \mu_m$. 
For an element $\bold Z_{(p)}/\bold Z$, we use the following notation 
$\< x \> = \frac{1}{e}\sum_{i=0}^{e-1}(<p^ix > -\frac{1}{2})$,
where $e$ is the minimal positive integer such that $p^e x =x$.
For $k=\bold F_q$,
we have the following proposition.
\proclaim{Proposition 8.2}
Let $\Cal G(D,\psi)$ be the hypergeometric sheaf of a hypergeometric data
$D$.  Let $\chi \in R \otimes I(\bar\bold F_q, \bar\bold Q_l)$
be a character of $\pi_1(\bold T(R))$ of finite order. We put
$$
C(\Cal G(D,\psi ), \chi ) =
\sum_{i=0}^{2(r-n-1)}(-1)^i
(\tr Fr_{\bold F_q} ;H_c^{i}(\bold T(R)\otimes \bar\bold F_q, 
\Cal G(D,\psi ) \otimes K(\chi )))
$$
Then we have $C(\Cal G(D, \psi ), \chi ) =
 \prod_{i=1}^r g(l_i(\chi ) + \kappa_i ,\psi )$.
\endproclaim
\demo{Definition}
For a constructible sheaf $\Cal F$ of $\bold T$,
the function from $\Hom (\bold T, \bold G_m) \otimes 
I(\bar\bold F_q, \bar\bold Q_l)$ to $\bar\bold Q_l$ given by
$$
\chi \mapsto C(\Cal F, \chi ) =
\sum_{i=0}^{2(r-n-1)}(-1)^i
\tr Fr_{\bold F_q} ;H_c^{i}(\bold T(R)\otimes \bar\bold F_q, 
\Cal F \otimes K(\chi ))
$$
is called the cohomological Mellin
transform of $\Cal F$.
\enddemo
\demo{Proof}
We can prove the proposition in the same way.  The variety $\Cal X$
$$
\Cal X =\{ (x,t,a) \mid \bold T(x) \times \bold A(t) \times \bold T(a)
\mid f(a,x) = t \}
$$
and the sheaf 
$R^n\varphi_*K(\Cal X;x^\alpha a^\kappa)\otimes \Cal L_{\psi}(t)$ on 
$\bold T(a)$
is descent to $\bold T(R)$.  Via the same change of coordinate, the 
descent variety $\Cal X_R$ is isomorphic to $\Cal X_b$;
$$
\Cal X_b =\{ (b, t ) \in \bold T(b) \times \bold A(t ) \mid
\sum_{i=1}^r b_i = t \}
$$
The descent sheaf on $\Cal X_b$ corresponds to the Kummer character
$\kappa (b^\kappa ) =$
\linebreak
$\kappa (\Cal X_b ; \prod_{i=1}^r b_i^{\kappa_i} )$.
Since $R^i \varphi_{!}K(\Cal X;x^\alpha a^\kappa)
\otimes \Cal L_{\psi}(t ) \simeq
R^i \varphi_{*}K(\Cal X; x^\alpha a^\kappa)
\otimes \Cal L_{\psi}(t)=0$
for $i \neq n$, we have
$H_c^{r-1+i}(\bold T(b), K(b^{\kappa + \bold l (\chi )})
\otimes \Cal L_{\psi}(t)) = 
H_c^{r-n-1+i}(\bold T(R), \Cal G(D, \psi) \otimes K(\chi ))$.
The alternating sum of the trace of Frobenius substitution on
the cohomologies of $\Cal X_b$ 
is given by
$$
\align
C(\Cal G(D,\psi ), \chi ) = &
\sum_{b_i \in \bold F_q^\times} \chi_{\kappa_i+l_i(\chi )}(b_i)
\psi(\sum_{i=1}^r\tr_{\bold F_q/\bold F_q}(b_i)) \\
= & 
\prod_{i=1}^r ( \sum_{b_i \in \bold F_q^\times} 
\chi_{\kappa_i+l_i(\chi )}(b_i)
\psi (\tr_{\bold F_q / \bold F_p}(b_i))) \\
= & \prod_{i=1}^r g(l_i(\chi ) + \kappa_i ,\psi ). \\
\endalign
$$ 
\enddemo
\proclaim{Corollary}
The order of $C(\Cal G(D,\psi ))$ is given by
$$
\ord_{\bar p} C(\Cal G(D, \psi )) = [\bold F_q:\bold F_p]
\sum_{i=1}^r ( \< l_i(\chi ) + \kappa_i \> +\frac{1}{2})
$$
\endproclaim
\heading
\S 9 A Lemma for independentness
\endheading
A hypergeometric data $D=(R, \{l_i\}_i \{\kappa\}_i)$ is called
primitive if $l_i(R)=\bold Z$.
Let $D = (R, \{l_i\}_i, \{\kappa_i\}_i)$ be a primitive hypergeometric
data.  The data $D$ is called non-divisorial if there exists no
$i,j$ ($i \neq j$) such that $l_i + l_j =0$ and $\kappa_i + \kappa_j =0$.
In this section, we prove the following proposition.
\proclaim{Proposition 9.1}
Let $D=(R, \{l_i\}_i, \{\kappa_i\}_i)$ and 
$D'=(R, \{l'_j\}_i, \{\kappa'_j\}_i)$ be primitive reduced non-divisorial
hypergeometric data and $c, c'$ be elements in $\bold Q$.
\roster
\item
Suppose that the equality
$$
\sum_{i=1}^r(<l_i(\chi)+\kappa_i>-\frac{1}{2})+c =
\sum_{i=1}^{r'}(<l'_i(\chi)+\kappa'_i>-\frac{1}{2}) +c'
$$
holds for a dense open set in $\chi \in R \otimes (\bold R/\bold Z)$.
(We say that the equality holds for almost all $r$ for short.)
Then $D$ is equal to $D'$ up to permutation and $c=c'$.
\item
Suppose that $\kappa_i, \kappa'_j \in I(\bar\bold F_q, \bar\bold Q_l)$ and
that the equality
$$
\sum_{i=1}^r\<\sigma ( l_i(\chi)+\kappa_i )\> + c =
\sum_{i=1}^{r'}\<\sigma ( l'_ i(\chi)+\kappa'_i) \> +c'
$$
holds for all $\chi \in R\otimes I(\bar\bold F_q, \bar\bold Q_l)$ and 
$\sigma :\bar\bold Q_l \to \bar\bold Q_l$.  Then $D$ is equal to
$D'$ up to permutation and $c=c'$.
\endroster
\endproclaim
To prove the proposition, we prove the following lemma.
\proclaim{Lemma 9.2}
\roster
\item
Let $\kappa_i \in \bold Q/ \bold Z$ be
distinct elements and $c \in \bold Q$.  If 
$$
\sum_{i=1}^r m_i (<x+\kappa_i>-\frac{1}{2}) + c =0
$$
holds for almost all $x \in \bold R /\bold Z$.  Then $m_i = 0$ for
all $i$ and $c=0$.
\item
Let $\kappa_i \in \bold Z_{(p)}/\bold Z$ be distinct elements
and $c\in \bold Q$. If 
$$
\sum_{i=1}^r m_i \< t(x+\kappa_i)\> +c = 0
$$ 
for all
$x \in \bold Z_{(p)}/\bold Z$ and $t \in (\hat\bold Z')^\times$,
then $m_i = 0$ for all $i$ and $c=0$.
Here $\hat\bold Z' =
\underset{\underset{(m,p)=1}\to{\longleftarrow}}\to{\lim}
\bold Z /m \bold Z$.
\endroster
\endproclaim
\demo{Proof}
(1) Let $f(x) = \sum_{i=1}^r m_i(<x+\kappa_i>-\frac{1}{2})+c$.
Since $m_i = \lim_{\epsilon \to +0}
f(x-\kappa_i+\epsilon) + f(x-\kappa_i-\epsilon)$, we get the theorem.
(2) To prove this theorem, we use an arithmetic argument.
We prove this statement in Appendix.
\enddemo
\demo{Proof of Proposition 9.1}
(1) We prove the proposition by the induction on rank $R$. 
If we take a base of $R \simeq \bold Z$, $l_i(\chi)$ and $l'_j(\chi)$ can
be written as $\chi$ or $-\chi$ by the primitivity condition.  Let $S_{\pm}$
and $S'_{\pm}$ be
$S_{\pm} = \{ i \mid l_i = \pm \chi\}$ and
$S'_{\pm} = \{ j \mid l'_j = \pm \chi\}$. By the assumption of 
non-divisoriality, $K_+=\{\kappa_i \mid i\in S_+ \}$ (resp.
$K'_+=\{\kappa_j \mid j\in S'_+ \}$ )
and $-K_-=\{-\kappa_i \mid i\in S_- \}$ 
(resp. $-K'_-=\{-\kappa_j \mid j\in S'_- \}$ ) does not intersect to each
other.  Since $<x>-\frac{1}{2} =-(<-x>-\frac{1}{2})$ for almost all 
$x$, we have
$$
\align
\sum_{i=1}^r(<l_i(\chi) + \kappa_i>-\frac{1}{2}) = &
\sum_{i \in S_+}(<\chi + \kappa_i>-\frac{1}{2}) -
\sum_{i \in S_-}(<\chi - \kappa_i>-\frac{1}{2}) \\
= &
\sum_{\kappa_i \in K_+}m_i (<\chi + \kappa_i>-\frac{1}{2}) -
\sum_{\kappa_i \in K_-}m_i (<\chi - \kappa_i>-\frac{1}{2}) \\
\endalign
$$
for almost all $\chi$, where $m_i= \# \{i \mid \kappa_i \in K_+ \cup K_-\}$
is the multiplicity of $\kappa_i$.
By the assumption,
$$
\align
& \sum_{\kappa_i \in K_+}m_i (<\chi + \kappa_i>-\frac{1}{2}) -
\sum_{\kappa_i \in K_-}m_i (<\chi - \kappa_i>-\frac{1}{2})  +c \\
= &
\sum_{\kappa_j \in K'_+}m'_j (<\chi + \kappa'_j>-\frac{1}{2}) -
\sum_{\kappa_j \in K'_-}m'_j (<\chi - \kappa'_j>-\frac{1}{2}) +c' \\
\endalign
$$
for almost all $\chi$.  Here $m'_j$ is the multiplicity of $\kappa'_j$
defined in the same way.  Since $m_i, m'_j >0$ for all $i,j$, we
have $c=c'$, $K_+ = K'_+$, $K_- = K'_-$ and
$m_i = m'_j$ if $\kappa_i = \kappa'_j$ by Lemma 9.2.

Assume that $\rank R = r -n-1 \geq 2$.  Since $l_1$ is primitive, we choose
a base $v_1, \dots , v_{r-n-1}$ such that 
$\Ker (l_1) = <v_2, \dots , v_{r-n-1}>$.  The coordinate is written as
$\chi_1, \dots ,\chi_{r-n-1}$.  Let $S=\{ i \mid l_i = \pm l_1\}$
and $S'=\{ i \mid l'_j = \pm l_1\}$.
For $i \notin S$ and $j \notin S'$, $\bar l_i$ and $\bar l'_j$
be the restriction of $l_i$ and $l'_j$ to $\Ker (l_1)$.  Since
$\bar l_i$ and $\bar l'_j$ is not zero, we have $\bar l_i = d_i \hat l_i$
and $\bar l'_j = d'_j \hat l'_j$ with primitive forms $\hat l_i$
and $\hat l'_j$.  Therefore for 
$r_1v_1 \in \bold Q v_1$, we have
$<\bar l_i+\bar l_i(\chi_1v_1) + \kappa_i>-\frac{1}{2} =
\sum_{k=1}^{d_i}
(<\hat l_i + \frac{1}{d_i}(\bar l_i (\chi_1v_1) + \kappa_i +k)>-\frac{1}{2})$.
Therefore $\sum_{i \notin S}(<\bar l_i + \bar l_i(\chi_1v_1) + \kappa_i>-
\frac{1}{2})$ can be expressed as a sum of $<L_l + K_l>$, where
$L_l$ is a primitive linear form on $\Ker (l_1)$.  
If necessary, by canceling terms,
we may assume that the sum does not contain divisorial pair 
in this expression. 
Since the restriction
of $\sum_{i \in S}(<l_i + \kappa_i> -\frac{1}{2})$ to
$\Ker (l_1) + \chi_1v_1$ is a constant function, we have
$$
\sum_{i \in S}(<l_i(\chi_1v_1) + \kappa_i> -\frac{1}{2}) +c =
\sum_{j \in S'}(<l'_j(\chi_1v_1) + \kappa'_j> -\frac{1}{2}) +c'
$$
for almost all $\chi_1$ by using the assumption of induction for $\Ker (l_1)$.
(Note that non-divisoriality condition inherits to the restriction to 
$\Ker (l_1) + \chi_1v_1$.)
Applying the assumption for $v_1 \bold Q$, $c=c'$,
$\{ \kappa_i \mid l_i = l_1 \} = \{ \kappa'_j \mid l'_j = l_1 \}$,
$\{ \kappa_i \mid l_i = -l_1 \} = \{ \kappa'_j \mid l'_j = -l_1 \}$ and
$\#\{ i \in S \mid \kappa_i = \kappa \} =
\#\{ j \in S' \mid \kappa'_j = \kappa \}$ for all 
$\kappa \in \bold Q/\bold Z$. Replacing $l_1$ to another primitive linear
form on $R$, we get the proposition.

For (2), we can prove the same argument.
\enddemo
\proclaim{Lemma 9.3}
Under the non-resonance condition, the primitive hypergeometric data
is non-divisorial
\endproclaim
\demo{Proof}
Assume that $l_i + l_j =0$ and $\kappa_i+ \kappa_j =0$.  
For simplicity, we may assume that $i =1, j=2$. Let
$R' = \Ker (l_1)$. Consider the following commutative diagram, where
$L' = \bold Z^{r-2}/R'$.
$$
\CD
@. 0 @. 0 @. 0 \\
@. @VVV @VVV @VVV \\
0 @>>> R' @>>> \bold Z^{r-2} @>>> L' @>>> 0 \\
@. @VVV @VVV @VVV \\
0 @>>> R @>>> \bold Z^{r} @>>> L @>>> 0 \\
@. @VVV @VVV @VVV \\
0 @>>> \bold Z  @>>> \bold Z^{2} @>>> \bold Z @>>> 0 \\
@. @VVV @VVV  \\
@. 0 @. 0  \\
\endCD
$$
Then we have $L/L' \simeq \bold Z$ and 
$$
L = L' \oplus \omega_1. \tag{*}
$$
Take an element $\rho \in R$ such that $l_1 (\rho) =1 $ and $l_2(\rho) =-1$.
Put $l_k (k) = a_k$. Then we have 
$$
\omega_1 -\omega_2 + \sum_{k=1}^r a_k \omega_k=0. \tag{**}
$$  The cone $C'$ generated 
by $\omega_i$ for $i=3,\dots ,r$ is contained in codimension 1 subspace
$L'\otimes R$ and the cone $C(\Delta )$ generated by 
$\omega_i$ ($i=1, \dots ,r$)
is a cone of $\omega_1$ and $\omega_2$ over $C'$ by the equality (**).
Therefore $C'$ is a codimension 1 face of $C(\Delta)$. It is easy to
see that the equation of for the face $C'$ is given by the projection
$pr$ to the $L'$ component via the decomposition (*).
Since the image $\alpha$ of $\kappa$ is given by
$$
\alpha =
\kappa_1\omega_1 + (-\kappa_1)\omega_2 + \sum_{i=3}^r \kappa_i\omega_i =
- \kappa_1\sum_{i=3}^ra_i\omega_i +\sum_{i=3}^r \kappa_i\omega_i,
$$
we have $pr (\alpha)=0$.  This contradicts to the non-resonance
condition.
\enddemo
\heading
\S 10 Main Theorem
\endheading
In this section we state and prove the main theorem.
First we prove the following lemma
\proclaim{Lemma 4.1}
Let $\Cal G(D)$ and $\Cal G(D, \psi)$ be hypergeometric sheaves on
$\bold T(R)$ for $k= \bold C$ and $\bold F_q$ respectively.
Let $\alpha$ be an element in $\bold T(R)(\bold C)$ 
(resp. $\bold T(R)(\bar \bold F_q)$)
such that $T_{\alpha}^*\Cal G(D) = \Cal G(D)$
(resp. $T_{\alpha}^* \Cal G(D, \psi ) = \Cal G(D, \psi )$) if
$k= \bold C$ (resp. $k=\bold F_q$).  Then we have $\alpha =1$.
\endproclaim
\demo{Proof}
First we assume $k= \bold C$.  The subgroup
$$
\Stab (\Cal G(D)) =
\{ \alpha \in \bold T(R)(\bold C) \mid T_{\alpha}^*\Cal G(D) \simeq 
\Cal G(D) \}
$$
is an algebraic subgroup because it is isomorphic to
$$
\Stab (\Cal G(D)) =
\{ \alpha \in \bold T(R)(\bold C) \mid T_{\alpha}^*DR(\Cal G(D)\mid_U) \simeq 
DR(\Cal G(D)\mid_U) \}.
$$
where $DR$ is a de Rham functor for local system on an open set $U$
of $\bold T$.  Therefore if $\Stab (\Cal G(D))$ is not zero, it contains
a non-trivial finite subgroup $S$.
Let $k=\bold F_q$. If 
$$
\Stab (\Cal G(D,\psi)) =
\{ \alpha \in \bold T(R)(\bar\bold F_q) \mid T_{\alpha}^*\Cal G(D,\psi) \simeq 
\Cal G(D, \psi) \}
$$
is non-trivial group, it contains a non-trivial finite subgroup $S$.
Therefore the sheaves $\Cal G(D)$ (resp.  $\Cal G(D, \psi)$) descend
to a sheaf $\Cal G(D)'$ (resp. $\Cal G(D, \psi)'$)
on $\bold T(R)/S$, since $\Cal G(D)$ and $\Cal G(D, \psi)$ are
irreducible.
Therefore we have
$$
\align
\chi_c (\Cal G(D)) & = \# (S)\cdot \chi_c(\Cal G(D)') \\
\chi_c (\Cal G(D,\psi )) & = \# (S)\cdot \chi_c(\Cal G(D,\psi)'). \\
\endalign
$$
This contradicts to the fact that
$\chi_c (\Cal G(D))= 1$ (resp. $\chi_c (\Cal G(D,\psi))= 1$).
\enddemo
Our main theorem is the following.
\proclaim{Theorem 4}
\roster
\item
Let $k=\bold C$.
Let $\Cal F$ and $\Cal F'$ be hypergeometric sheaves on
a torus $\bold T$.
If there exists an isomorphism
$\phi : \Cal F\mid_U \to \Cal F'\mid_U$ 
as a variation of $K$-Hodge structures on an open set $U$ in $\bold T$,
then there exists an algebraic correspondence $Al$ and a
Hodge cycle of Fermat motif $F$ such that $\phi = Al\cdot F$.
\item
Let $k=\bold F_q$.
Let $\Cal F$ and $\Cal F'$ be hypergeometric sheaves on
a torus $\bold T$.
If there exists an isomorphism
$\phi : \Cal F\mid_U \to \Cal F'\mid_U$ 
as a $\bold Q_l$ local system on an open set $U$ in $\bold T$,
then there exists an algebraic correspondence $Al$ and a
Tate cycle of Fermat-Artin-Shreier motif $FAS$ such that 
$\phi = Al\cdot FAS$.
\endroster
\endproclaim
(1) 
Let $D$ be a hypergeometric data.
By using multiplicative correspondence successively, there exist
a primitive hypergeometric data $D^{(p)}$ and elements 
$\lambda_1, \dots, \lambda_a$ such that 
there exists an isomorphism 
$\Cal G(\tilde D) \simeq T^*_\alpha \Cal G(D^{(p)})$ 
for some $\alpha \in\bold T(R)$
induced by an algebraic 
correspondence, where $\tilde D = D \oplus \oplus_{i=1}^a(0,\lambda_i)$.
Let $D^{(p)}_{red}$ be the reduced part of $D^{(p)}$. Then 
there exist Fermat motives $F_1, F_2$ such that 
$\Cal G(D) \otimes F_1 \simeq T^*_\alpha\Cal G(D^{(p)}_{red}) \otimes F_2$
for some $\alpha \in \bold T(R)$.
Therefore there exist a Fermat motif  $F_3$ and
the following isomorphism induced by 
the composite of algebraic correspondences
$$
\Cal G(D) \otimes F_3 \simeq T_{\alpha}^*\Cal G(D^{(p)}),\tag{*}
$$
for some  $\alpha \in bold T(R)$.
By the isomorphism (*), the difference of Mellin transform of $\Cal G(D)$
and $\Cal G(D^{(p)}_{red})$ is a constant.
Now we compare cohomological Mellin transform of sheaves $\Cal F$ and
$\Cal F'$ on $\bold T=\bold T(R')$.  Since $\Cal F$, $\Cal F'$ 
$\Cal F\mid_U$ and $\Cal F'\mid_U$
are cohomological
mixed Hodge complices, there is natural mixed Hodge structure on
\linebreak
$\Hom (\Cal F, \Cal F') \simeq \bold Q$
and $\Hom (\Cal F\mid_U, \Cal F'\mid_U) \simeq \bold Q$.
Therefore the restriction map 
$$
\Hom (\Cal G(D), \Cal G(D')) \to \Hom (\Cal G(D)\mid_U, \Cal G(D')\mid_U)
$$
is an isomorphism of Hodge structure.  Since the right hand
side is pure of weigh $0$ with type $(0,0)$, $\Cal F$ and $\Cal F'$
are isomorphic as a mixed Hodge complex.
Therefore the Mellin transform 
$C(\Cal F, \chi )$ and $C(\Cal F', \chi)$ of 
$\Cal F$ and $\Cal F'$ coincide for all 
$\chi \in R'' = \Hom( \bold T, \bold C^{\times})\otimes \bold Q$.
By the definition of hypergeometric sheaves, there exists an isomorphism
$\Phi : \bold T \simeq \bold T(R)$ and $\Phi' :\bold T \simeq \bold T(R')$ 
and hypergeometric data $D=(R, \{l_i\}_i, \{\kappa_i\}_i)$ and
$D'=(R', \{l'_i\}_i, \{\kappa'_i\}_i)$ such that
that $\Cal F = \Phi^* \Cal G(D)$ and $\Cal F' = {\Phi'}^*\Cal G(D')$.
Via the isomorphism $\Phi$ and $\Phi'$, we have an isomorphism
of $R \simeq R' \simeq R''$. Under this isomorphism, we identify
$R$, $R'$ and $R''$.
Via this identification, cohomological Mellin transform $\Cal G(D)$ and
$\Cal G(D')$ coincide with those of $\Cal F$ and $\Cal F'$.
Therefore the Mellin transform of $\Cal G(D^{(p)}_{red})$ and
$\Cal G({D'}^{(p)}_{red})$ coincide up to constant.
By using Proposition 8.1,
the explicit calculation of the Hodge type
of Mellin transform of $\Cal G(D^{(p)}_{red})$ and
$\Cal G({D'}^{(p)}_{red})$, we have
$$
\align
& (\sum_{i=1}^r<\sigma(\kappa_i + l_i(\chi ))> + 
<-\sum_{i=1}^r \sigma\kappa_i> \\
= &
(\sum_{i=1}^r<\sigma(\kappa'_i + l'_i(\chi ))> + 
<-\sum_{i=1}^r \sigma\kappa'_i> +c, \\
\endalign
$$
for some constant $c$,
where $D^{(p)}_{red}=(R, \{l_i\},\{\kappa_i\})$ and
${D'}^{(p)}_{red}=(R, \{l'_i\},\{\kappa'_i\})$.
Therefore by Proposition 9.1 , hypergeometric data $D^{(p)}_{red}$ is equal to
${D'}^{(p)}_{red}$ up to permutation.
Again using (*), there exists a Fermat motif $F_4$ and an isomorphism
$$
\Cal G(D) \otimes F_1 \simeq T_\alpha^*\Cal G(D')
$$
induced by an algebraic correspondence.  Comparing to the assumption,
$F_1$ is a Hodge cycle.  Since $\Hom (\Cal G(D), \Cal G(D'))$ is
one dimensional, this is equal to the original isomorphism up to 
constant. Since $T^*_{\alpha}\Cal G(D) = \Cal G(D)$,
implies $\alpha = 1$, we have the theorem.

(2) The proof is similar. 
The isomorphism $\Cal G(D^{(p)}_{red},\psi )\mid_U
\simeq \Cal G({D'}^{(p)}_{red},\psi )\mid_U $
implies
$$
(\sum_{i=1}^r \< t(\kappa_i + l_i(\chi )) \> + 
\< -\sum_{i=1}^r t\kappa_i \>
=
(\sum_{i=1}^r \< t(\kappa'_i + l'_i(\chi )) \> + 
\< -\sum_{i=1}^r t\kappa'_i \>,
$$
for all $t\in (\hat\bold Z')^\times$ and 
$\chi \in R \otimes (\bold Z_{(p)}/\bold Z)$
by taking $\bar p$-adic order of Mellin transform (Corollary to 
Proposition 8.2).  This implies $D^{(p)}_{red}$ is equal to
${D'}^{(p)}_{red}$ up to permutation by Proposition 9.1 (2).
The rest of the proof is similar to that of (1).

\heading
\S Appendix  The proof of Proposition 9.1 (2)
\endheading
The following lemma is well known.
\proclaim{Lemma A.1}
Let $N$ be a natural number greater than 2, $\psi$ be a additive character
of $\bold Z/N\bold Z$ given by $\psi (x) = \exp (2\pi\sqrt{-1}x/N)$.
Then for an odd multiplicative character $\chi$, 
$$
L_N(\chi , 1) = -\frac{2\pi\sqrt{-1}}{N}
\sum_{m,l=0}^{N-1}
\psi (ml)\chi (m)(<\frac{l}{N}>-\frac{1}{2})
$$
is non-zero.
\endproclaim
\demo{Proof}
See [BS].
\enddemo
From now on we assume that $p \neq 2$ and $N$ is even.
We introduce a set $S$ (super singular divisor of $N$ for $p$)
and $S^c$
of divisor of $N$ by
$$
\align
S=\{ M \mid & M \text{ divides } N, \text{ the subgroup $<p>^{\times}$ of 
$(\bold Z/M\bold Z)^{\times}$} \\
& \text{ generated by $p$ contains $-1$ },
M \neq 1,2
\}. \\
S^c=\{ M \mid &  M \text{ divides } N, \text{ the subgroup $<p>^{\times}$ of 
$(\bold Z/M\bold Z)^{\times}$ } \\
& \text{ generated by $p$ does not contain $-1$ },
M \neq 1,2
\}. \\
\endalign
$$
By definition of $\< \>$, we can check the following property.
\roster
\item $\< 0\> = - \frac{1}{2}$
\item
If $x \neq 0$,
$\< -x \> = - \< x \>$.
\item
If $x \in \frac{1}{M}\bold Z/\bold Z - 0$ with $M \in S$, then
$\< x \> =0$.
\endroster
The properties (1) and (2) is direct consequence from the definition.
To show (3), use the equality $\< px \> = \< x \>$ and (2).
Let $Func(S)$ be the vector space of $\bold C$-valued function on the
set $S$.
The natural linear map from $Func(\frac{1}{N}\bold Z/\bold Z)$ to
$Func((\bold Z/N\bold Z)^{\times} \times (\frac{1}{N}\bold Z/\bold Z))$
is induced by the second projection.  Let $V$ be the quotient space
$$
Func((\bold Z/N\bold Z)^{\times} \times (\frac{1}{N}\bold Z/\bold Z))
/Func(\frac{1}{N}\bold Z/\bold Z).
$$
We define a linear map $\phi$
from $\oplus_{\kappa \in \frac{1}{N}\bold Z/\bold Z}
\bold C v_{\kappa}$ to $V$ defined by 
$\phi (v_{\kappa}) = \< t(x + \kappa)\>$.
Here $\< t( x+ \kappa )\>$ is the class of function
on $(t, x) \in (\bold Z /N\bold Z)^{\times} \times 
(\frac{1}{N}\bold Z /\bold Z)$
in $V$.
For $\kappa \in \frac{1}{N}\bold Z/\frac{M}{N}\bold Z, M \in S\cup\{2\}$
we define an element $\sigma_{M,\kappa}$ by
$$
\sigma_{M,\kappa}=\sum_{i=0}^{l-1}v_{\kappa + \frac{i}{l}}.
$$
where $l = N/M$.  Then $\sigma_{M,\kappa } \in \Ker (\phi )$.
In fact,
$$
\align
\phi (\sigma_{M,\kappa }) = &
\sum_{i=0}^{l-1}\< t(x+\kappa + \frac{i}{l})\> \\
= &
\< lt(x+\kappa ) \>.  \\
\endalign
$$
Since $lt(x+\kappa) \in \frac{1}{M}\bold Z /\bold Z$,
it is 0 if $x \neq -\kappa \text{ mod }\frac{1}{l}\bold Z$ 
and $-\frac{1}{2}$ if $x = -\kappa \text{ mod }\frac{1}{l}\bold Z$.
Therefore the class of this function in $V$ is 0. 
Now we choose a numbering of 
$S\cup\{ 2\}=\{ M_1, \dots , M_a\}$ such that the
the divisibility $M_i \mid M_j$ implies $j \leq i$.
Let $t_j = \sum_{i=1}^j\varphi (M_i)$.  
\proclaim{Proposition A.2}
\roster
\item
$\dim \Im (\phi) \geq \sum_{M \in S^c}\varphi (M)$.
\item 
Let $W$ be the
subspace of $\Ker (\phi )$ generated by $\sigma_{M_i,t_{i-1}+s}$
($i=1, \dots , a, s=1, \dots ,\varphi (M_i)$), 
$\sigma_{2,0}$ and $\sigma_{2,1}$.
Then
$\dim W  \geq \sum_{M \in S}\varphi (M) +2$.
\endroster
As a consequence, the equalities holds for both (1) and (2)
and $\sigma_{M,\kappa }$
($\kappa \in \frac{1}{N}\bold Z/\frac{1}{l}\bold Z,$ 
\linebreak
$M \in S$),
$\sigma_{2,0}$ and $\sigma_{2,1}$
generates $\Ker (\psi )$.
\endproclaim
\demo{Proof}
(1).
For an element $M \in S^c$, there exists an odd character $\chi$ of
$(\bold Z /M \bold Z)^{\times}$ such that $\chi (p) =1$.
For an element $l' \in \bold Z/N\bold Z$ such that
$l'(\frac{1}{N}\bold Z/\bold Z) = \frac{1}{M}\bold Z/\bold Z$,
we define a linear map $\theta_{l'}$ from $V$ to $\bold C$ by
$$
\theta_{l'}(f) = N_{M,l'} \sum_{t \in (\bold Z/N \bold Z)^{\times},
x \in \frac{1}{N}\bold Z/\bold Z}\chi (t) \psi (l'x) f(t,x),
$$
where $N_{M, l'}$ is the normalizing factor;
$$
N_{M, l'} = \chi (l/l') \frac{M\varphi (M)}{2\pi\sqrt{-1}\varphi (N)}
L_M(\chi ,1)^{-1}.
$$
Then $\theta_{l'}\circ \phi (v_{\kappa})$ is calculated as follows.
$$
\align
\theta_{l'}\circ \phi (v_{\kappa}) = &
N_{M,l'} \sum_{t \in (\bold Z/N \bold Z)^{\times},
x \in \frac{1}{N}\bold Z/\bold Z}\chi (t) \psi (l'x) \< t(x+\kappa \> \\
= &
N_{M,l'} \psi (-l'\kappa )
\sum_{t \in (\bold Z/N \bold Z)^{\times},
x \in \frac{1}{N}\bold Z/\bold Z}\chi (t) \psi (l'x) \< tx \>  \\
= &
N_{M,l'} \psi (-l'\kappa )
\sum_{t \in (\bold Z/N \bold Z)^{\times},
x \in \frac{1}{N}\bold Z/\frac{1}{l}\bold Z}
\chi (t) \psi (l'x) (\sum_{y \text{ mod }\frac{1}{l}\bold Z/\bold Z = x}
\< ty \>)  \\
= &
N_{M,l'} \psi (-l'\kappa )
\sum_{t \in (\bold Z/N \bold Z)^{\times},
x \in \frac{1}{N}\bold Z/\frac{1}{l}\bold Z}
\chi (t) \psi (l'x) \< tlx \>  \\
\endalign
$$
By changing summation $t = t_0 \tau, t_0 x = \xi$
where $t_0l = l'$, it is equal to
$$
\align
& 
N_{M,l'} \psi (-l'\kappa )
\sum_{\tau \in (\bold Z/N \bold Z)^{\times},
l \xi \in \frac{1}{M}\bold Z/\bold Z}
\chi (t_0 \tau ) \psi (l \xi ) \< \tau l \xi \>  \\
= &
\psi (l'\kappa ) N_{M,l'}
\chi (l'/l) \frac{2\pi\sqrt{-1}\varphi (N)}{M\varphi (M)}
L_M(\chi ,1) \\
= &
\psi (- l'\kappa ).  \\
\endalign
$$
For $M \in S^c$, there exist $\varphi(M)$ elements $l'$ in
$\bold Z/N\bold Z$ satisfying (*).
By considering $\theta_{l'}$, for all $M \in S^c$ and
$l'$ with the condition, we get $\sum_{M \in S^c}\varphi (M)$
linear homomorphism which is independent on the image of $\psi$.
In fact the determinant of
$(\theta_{l'}(v_{\kappa}))_{(M,l'),\kappa } = 
(\psi (-l'\kappa))_{(M,l'),\kappa }$
does not vanish by Van der Monde determinant.

(2) We have already shown that $\sigma_{M,\kappa }$ is contained in
$\Ker (\phi )$.
Let $M \in S$ and $l'\in \bold Z/N\bold Z$ such that
$l' (\frac{1}{N}\bold Z/\bold Z) = \frac{1}{M}\bold Z/\bold Z$.
To show that there exist enough independent $\sigma_{M,\kappa }$,
we define linear forms $\gamma_{M,l'}$ for $(M,l')$ satisfying
the above condition.
Since $-1 \neq 1 \in (\bold Z/M\bold Z)^{\times}$, we choose an
odd character $\chi $ of $(\bold Z/M\bold Z)^{\times}$
Let us define $\gamma_{M,l'}$ by
$$
\gamma_{M,l'}(v_\kappa) = N_{M,l'}\sum_{t\in (\bold Z/M\bold Z)^{\times}, 
x \in \frac{1}{N}\bold Z/\bold Z}\chi (t) \psi (l'x)
(<t(x+\kappa )>-\frac{1}{2}),
$$
where the normalizing factor $N_{M,l'}$ is given by 
the same formula as before. Then we can show that
$$
\gamma_{M,l'}(v_{\kappa})= \psi (-l'\kappa).
$$
Let $\tilde M \in S$ and $\tilde l \tilde M =N$. Then we have
$$
\align
\gamma_{M,l'}(\sigma_{\tilde M,\tilde \kappa}) = &
\sum_{i=1}^{\tilde l-1}\psi (-l'(\kappa + \frac{\tilde M}{N}i)) \\
= & \cases 
\tilde l \psi (-l'\kappa ) & \qquad\text{ if }M \mid \tilde M \\
0 & \qquad\text{ otherwise } \\
\endcases
\endalign
$$
Then for $M_j \in S$, the determinant
of $(\gamma_{M_j,l'}(\sigma_{M_j,t_{j-1}+s}))_{l',s=1, \dots ,\varphi (M_j)}$
is equal to 
$$
\prod_{l'}\psi( -l'(t_{i-1}+1))
\prod_{l'<l''}(\psi(-l')-\psi (-l''))
$$
and non-zero. Since 
$$
(\gamma_{M_i,l'}(\sigma_{M_j,t_{j-1}+s}))_{l',s=1, \dots ,\varphi (M_j)}
$$
is a zero matrix if $i<j$.  Therefore the matrix
$$
(\gamma_{M,l'}(\sigma_{M, i}))_{(M,l),i=1, \dots ,\sum_{M\in S}\varphi (M)}
$$
is a non singular matrix.  Moreover, we can see $\sigma_{2,0}$
and $\sigma_{2,1}$ is annihilated by all $\gamma_{M,l'}$ ($M \in S$)
and independent.  Therefore we have prove the proposition.
last statement is the consequence of (1) and (2).
\enddemo  
Now we take a sufficiently large $N$ to prove the independentness of
$\< t(x+\kappa)\>$.
\demo{Proof of Proposition 9.1 (2)}
Let $N$ be an even integer such that 
$\kappa_i \in \frac{1}{N}\bold Z/\bold Z$ ($i=1, \dots ,r$). 
We choose a numbering of $S=\{M_1, \dots ,M_a\}$
satisfying the condition given just before Proposition 9.2.
Suppose that
$\sum_{i=1}m_i\< t( x+\kappa_i)\> + c = 0$.
Since the equality holds as a function of 
$t \in (\bold Z/N\bold Z)^{\times}, x \in \frac{1}{N}\bold Z/\bold Z$,
modulo constant, $\sum_{i=1}^rm_iv_{\kappa_i}$
is contained in the submodule generated by $\sigma_{M_i, j}$
($i= 1, \dots , a, j=t_{i-1}+1, \dots ,t_i$)
and $\sigma_{2,1}$ and $\sum_{i=0}^{N-1}v_i$
by Proposition A.2.  Suppose that
$\sum_{i=1}^rm_iv_{\kappa_i} = \sum a_{M_i,i}\sigma_{M_i,j}$
and 
$\sum a_{M,\kappa}(\sum_{k=0}^{M/N-1}\< t( x+ j + \frac{M_i}{N}k)\>)=0$
for $\frac{1}{CN}\bold Z/\bold Z$ for sufficiently large $C$.
Let $M_1 \geq 2$ be a maximal element where $a_{M_i,j} \neq 0$.
We choose $C$ such that
the subgroup generated by $p$ in $(\bold Z/CM_1\bold Z)^{\times}$
does not contain $-1$.  We apply the linear map $\theta_{l'}$
for $l'$ such that 
$l'(\frac{1}{CN}\bold Z/\bold Z) = \frac{1}{CM_1}$
by choosing an odd character $(\bold Z/NC\bold Z)^{\times}$
vanishing on the subgroup generated by $p$.
Then $\theta_{l'}(\sum a_{M_i,j}\sigma_{M_i,j})=
\sum_{\kappa}a_{M_i,j}\psi_{CN} (-l'j )$.
Here we used $\psi_{NC} (x) =\exp (2\pi \sqrt{-1}x/NC)$.
Since the determinant $(\psi_{CN} (-l'j))_{l'j=t_{i-1}+1,\dots ,t_i}$
is not zero, $a_{M_i,j}$ 
($j=t_{i-1}+1, \dots ,t_i$)
must be $0$. 
Suppose that 
$a\sum_{i \in \frac{1}{N}\bold Z/\bold Z}
\< t(x+i) \> +c=0$.  This implies $a\< Ntx \> +c = 0$ for all
$t$ and $x$.  By putting $x=\frac{1}{2N}$(resp. $x=0$), 
we have $c=0$ (resp. $-\frac{1}{2} a +c=0$). Therefore $a= c= 0$.
This is contradiction. 
\enddemo

\enddocument